\documentclass[11pt]{article}
\usepackage{amsmath,latexsym,amssymb,amsfonts,amsbsy, amsthm}
\usepackage [latin1]{inputenc}
\usepackage{color, xcolor}

\usepackage{graphicx}
\usepackage{soul, color}
\usepackage[margin = 1.5in]{geometry}

\newtheorem{theorem}{Theorem}[section]

\newtheorem{lemma}[theorem]{Lemma}

\newtheorem{remark}[theorem]{Remark}

\newcommand{\p}{\ensuremath{\partial}}

\newcommand{\eps}{\ensuremath{\varepsilon}}

\newcommand{\bigO}{\mathcal{O}}

\newcommand{\la}{\lambda}

\def\mau{\mathbf{u}}

\def\dv{\mathrm{div}}
\let\f=\frac

\let\va=\varepsilon

\let\f=\frac

\def\dv{\mathrm{div}}

\begin{document}

\title{\textbf{Structure stability of steady supersonic shear flow with inflow boundary conditions}}

\author{ {\bf Song $\text{Jiang}^{b,c}$ \qquad\quad Chunhui $\text{Zhou}^a$}
\\[4mm]
\small $^a$ Department of Mathematics, Southeast University,
Nanjing, 210096, China.\\
\small Email: zhouchunhui@seu.edu.cn\\
\small $^b$ Institute of Applied Physics and Computational Mathematics,\qquad\quad\quad\\
\small  Beijing 100088, China. Email: jiang@iapcm.ac.cn\qquad\quad\quad\\
\small $^c$ National Key Laboratory of Computational Physics, Beijing 100088, China.}
\date{}
\maketitle

\begin{abstract}
We study the existence and zero viscous limit of    smooth solutions to  steady compressible Navier-Stokes equations near plane shear
flows   between two moving parallel walls.
Under the assumption  $0<L\ll1$,  we prove that for any plane supersonic shear flow $\mathbf{U}^0=(\mu(x_2),0)$, there exist  smooth solutions near  $\mathbf{U}^0$ to  steady compressible Navier-Stokes equations in a 2-dimension domain $\Omega=(0,L)\times (0,2)$. Moreover, based on the uniform-in-$\va$ estimates, we establish the zero viscosity limit of  the solutions obtained above
to the solutions of the steady Euler equations.
\end{abstract}

\bigskip
\renewcommand{\theequation}{\thesection.\arabic{equation}}
\setcounter{equation}{0}

\section{Introduction}

In this paper, we shall study the structure stability of   steady supersonic  shear flows with inflow boundary condition in a 2-dimension domain $\Omega=(0,L)\times (0,2)$.  We are interested in the  steady compressible Navier-Stokes equations   in the following dimensionless form:
\begin{eqnarray}
\dv(\rho^\va\mau^\va)=0&&\quad\text{in}~\Omega,\label{0.3}\\
\rho^\va\mau^\va\cdot\nabla u^\va- \mu\va\Delta u^\va-\lambda\va\partial_{x_1}\dv\mau^\va+ \partial_{x_1}P^\va=0&&\quad\text{on}~\Omega,\label{0.4}\\
\rho^\va\mau^\va\cdot\nabla v^\va- \mu\va\Delta v^\va-\lambda\va\partial_{x_2}\dv\mau^\va+ \partial_{x_2}P^\va=0&&\quad\text{on}~\Omega.\label{0.5}
\end{eqnarray}
here $\va=1/Re,\ Re$ is the Reynolds number; $\mau^\va, \ \rho^\va$ are the velocity and the density, $P^\va$ is the pressure for isentropic flows  given by
$P^\va(\rho^\va)=a(\rho^\va)^\gamma$ with $a$ being a positive constant and
$\gamma>1$ being the specific heat ratio,  $\mu>0,\mu'>0$ are the scaled shear and bulk viscosity, $\lambda=\mu'+\f13\mu$, $c=\sqrt{(P^\va)'}$  the speed of sound. Without loss of generality, we will assume $\mu\equiv1$ in the following.  $\partial\Omega$ is divided into the inflow part
$\Gamma_\mathrm{in}$ ($\mathbf{u}\cdot\mathbf{n}<0$), the outflow part
$\Gamma_\mathrm{out}$ ($\mathbf{u}\cdot \mathbf{n}>0$), and the impermeable wall
$\Gamma_0$ and $\Gamma_2$ ($\mathbf{u}\cdot \mathbf{n}=0$). More precisely,
$$\Gamma_\mathrm{in}=\{x_1=0,\;\; 0\leq x_2\leq 2\}$$
$$\Gamma_\mathrm{out}=\{x_1=L,\;\; 0\leq x_2\leq 2\}$$
$$\Gamma_0=\{0\leq x_1\leq L,\;\; x_2=0\},\ \Gamma_2=\{0 \leq x_1\leq L,\;\; x_2=2\}.$$
The Navier-Stokes equations for a steady isentropic compressible viscous
flow is a mixed system of hyperbolic-elliptic type, as the momentum
equations are an elliptic system in the velocity, while the continuity
equation is hyperbolic in the density. Therefore, if we consider the
inflow boundary problem, it is necessary to prescribe the density on
the part of inflow boundary ($\mathbf{u}\cdot\mathbf{n}<0$). Besides,
We will consider the non-slip boundary condition on the moving walls. The boundary conditions under consideration are:
 \begin{eqnarray}
&& \mathbf{u}=\mau_{in},\ \text{on}~\Gamma_\mathrm{in},\ \mathbf{u}=\mau_{out},\ \text{on}\ \Gamma_\mathrm{out}, \
\mau\equiv(V_0, 0),\ \text{on}~\Gamma_0, \ \mau\equiv (V_1,0),\ \text{on}~\Gamma_2.\nonumber\\
&&\rho=\rho_0,\ \text{on}~\Gamma_\mathrm{in}.\label{b}
 \end{eqnarray}
Here $\mau_{in}(x_2)=(u_{in},v_{in})$,  $\mau_{out}(x_2)=(u_{out},v_{out})$, $V_0,V_1$ are positive constants. We also assume the following compatibility conditions on the corners:
\begin{eqnarray}
&&u_{in}(0)=u_{out}(0)=V_0, \ u_{in}(2)=u_{out}(2)=V_1,\nonumber\\
&& v_{in}(0)=v_{in}(2)=v_{out}(0)=v_{out}(2)=0.\label{com}
\end{eqnarray}

 The study on the asymptotic behavior of solutions of Navier-Stokes
equations as $\va\rightarrow0^+$  has been one of the most fundamental problems in fluid dynamics. In the steady setting of incompressible flow, if there is a mismatch between the basic Euler flows and the non-slip boundary conditions on the boundary,  there would be a thin fluid boundary layer of
size $\va^{\f12}$ to connect the Euler velocity profiles and the non-slip boundary conditions. The authors in  \cite{GM,GN,GI,GI1,I,IM,IM1} have established    the validity of the Prandtl boundary layer expansion and its error estimates. On the other hand, if there is no mismatch between the basic Euler flows and the non-slip boundary conditions on the boundary, then there would be no strong boundary layers near the rigid walls. The class of strictly parallel flows  satisfying the steady incompressible Navier-Stokes equations are limited, this includes two important special cases: the plane Couette  flow and the
 plane Poiseuille flow. The authors in \cite{IZ,JZ} have proved the existence and zero viscous limit of solutions near shear flows of Poiseuille-Couette type with non-slip boundary conditions on the rigid walls in Sobolev space. For the dynamic stability  of incompressible shear  flows or the boundary layer type flows with large Reynolds numbers we refer to \cite{BGM,BGM1,CLWZ,CWZ,GMM,GGN1,GGN,MZ,LMZ,R,WZ} and the reference therein.

  Considering the dynamic stability of compressible flows  with large Reynolds numbers, the authors in \cite{ADM,LZ,YZ1,YZ,WWZ,WW,ZZ} and among others have studied the linear stability of Poiseuille-Couette type flows or Prandtl type flows under different boundary conditions. In the steady setting, the authors in \cite{LZH} proved the existence and zero viscous limit of plane steady slightly compressible Navier-Stokes equations with Navier-slip boundary conditions. Recently, the authors in \cite{GW} studied the steady prandtl expansion for full
 compressible Navier-Stokes system  with non-slip boundary condition on the rigid wall and viscous-inflow boundary conditions in the flow direction under the assumption that the  Mach number is small. And the authors in \cite{LYZ} studied the structural stability of boundary layers in the entire subsonic regime in 2-D with non-slip boundary condition on the rigid wall and periodic conditions on the flow direction.

 For any parallel flow $\bold{U}^0=(\mu(x_2),0)$,
If we take $P^0\equiv C$ , here $C>0$ is a constant, then it is easy to check that $(P^0,\bold{U}^0)$
satisfy the stationary incompressible Euler equations:
\begin{equation} \label{main.euler}
\begin{cases}
\bold{U}^0 \cdot \nabla \bold{U}^0 + \nabla P^0 = 0 \\
\nabla \cdot \bold{U}^0 = 0.
\end{cases}
\end{equation}
If we consider the  compressible perturbation of viscous fluids near $(P^0,\bold{U}^0)$ in the absence of external forces, then we can obtain steady compressible Navier-Stokes equations (\ref{0.3})-(\ref{0.5}).
In this situation any parallel flows other than the plane  Couette flow are not solutions to system (\ref{0.3})-(\ref{0.5}). However, there are a large number of cases where the flow is essentially parallel to one direction, e.g., the inlet flow between parallel walls and flow along a flat plate.

   We are interested in the existence of solutions to steady compressible viscous flows around supersonic shear flows as well as the zero viscosity limit from steady compressible Navier-Stokes equations to steady incompressible Euler equations. For this purpose, we will first   expand the solution in $\va$ as:
\begin{equation}\label{app.0} \begin{cases}
u^\eps=\mu+\va u_e^1+\va u_p^1+\va^{\f32} u_e^2+\va^{\f32} u_p^2 + u \triangleq u_s +   u\\
v^\eps= \va v_e^1+\va^{\f32}v_p^1+\va^{\f32} v_e^2+\va^{2}v_p^2+ v\triangleq v_s +   v,\\
\rho^\eps =\rho^*+\va \rho_e^1+\va\rho_p^1+ \va^{\f32} \rho_e^2+\va^{\f32}\rho_p^2+\rho\triangleq \rho_s +   \rho,\\
\end{cases}
\end{equation}
Here $(u_e^i,v_e^i,\rho_e^i)$ and $(u_p^i,v_p^i,\rho_p^i)$, $i=1,2$, are defined in section 2.    We denote by
 \begin{equation}p^0=\min(\f83,q)>2,\label{p0}\end{equation}
here $q$ is defined as following:
$$q=\sup\{P^*|\ |2-2/P^*|<\la_*(\f\pi2)\}, \ \la_*(\f\pi2)>1.$$
For more details about $P^*,\la_*$ one can refer to Chapter 3 in \cite{K-M-R}. The main result of this paper reads as:
\begin{theorem}\label{main}
For $2<p<p^0$,  $\mathbf{U}^0=(\mu(x_2),0)$, $\rho^*>0$,   we assume that $\mu(x_2)\in C^6([0,2])$ satisfying:
\begin{equation}
\mu(0)=V_0>0,\ \mu(2)=V_1>0,\ \mu^2>c^2=(P^\va)'(\rho^*),\label{mu}
\end{equation}
and there is no mismatch between the basic flow  $\mathbf{U}^0$ and the moving boundaries, then there exists a triple  $(u_s,v_s,\rho_s)$ defined in (\ref{app.0}) such that if
 \begin{eqnarray}|
 &&\mau_{in}-\mau_s(0,\cdot)|_{C^2([0,2])}+|\mau_{out}-\mau_s(L,\cdot)|_{C^2([0,2])}+|\rho^0-\rho_s(0,\cdot)|_{C^2([0,2])}
 \nonumber\\
 \leq &&C \va^{\f52-\f2p+\sigma},\label{0.13}\end{eqnarray}
there exists a unique solution $(\mau^\va,\rho^\va)\in W^{2,p}(\Omega)\times W^{1,p}(\Omega)$ to the system (\ref{0.3})-(\ref{com}) with the remainder solution
$(u,v,\rho)$ defined in (\ref{app.0}) satisfying the following estimates:
\begin{eqnarray}
&&\|\mau\|_{L^2}+\va^{\f12}\|\nabla\mau\|_{L^2}+\|\rho\|_{L^2}\leq C\va^{\f52-\f2p+\sigma},\label{th.0.1}\\
&&\va\|\nabla^2\mau\|_{L^p}+\|\rho\|_{W^{1,p}}\leq C\va^{1+\f\sigma2}.
\end{eqnarray}
Consequently we have
\begin{eqnarray}
&&\|u^\va-\mu\|_{L^\infty}+\|v^\va\|_{L^\infty}+\|\rho^\va-\rho^*\|_{L^\infty}\leq C\va,\label{th.0.4}\\
&&\|\nabla\mau^\va-\nabla\mathbf{U}^0\|_{L^\infty}\leq C\va^{\f\sigma2},
\end{eqnarray}
where $\sigma>0$ is a constant sufficiently small, and the constant $C$ does not depend on $\va$.
\end{theorem}
\begin{remark}
The basic flow $\mathbf{U}^0=(\mu(x_2),0)$ in Theorem \ref{main}  is any shear flow satisfying the supersonic assumption (\ref{mu}).
\end{remark}

\begin{remark}
The constant $p^0$ in Theorem \ref{main} is not optimal. In fact, if we expand more terms in section 2, then $p^0$ can be taken to be any constant satisfying $2<p^0<q$.
\end{remark}

Let us make a few comments on the proof of Theorem  \ref{main}.
 The basic flows $\bold{U}^0 = (\mu(x_2), 0)$  in this paper satisfy the no-slip boundary conditions:
$\mu(0) = V_0,\ \mu(2)=V_1$ on the moving boundaries and there is no mismatch between the basic flows and the moving boundaries. Thus, there would be no strong boundary layers around  $y=0$ and $y=2$. However, weak boundary layers will still arise due to viscous effects. To obtain the existence of solutions to system (\ref{0.3})-(\ref{0.5}), we will use the multi-scale expansion of $(\rho^\va,u^\va,v^\va)$ around the basic flow $\bold{U}^0$.
First,  compressible Euler correctors
$(u_e^i,v_e^i,\rho_e^i)$ with slip boundary conditions are constructed to balance the shear stress perturbations. Under the supersonic assumption of the basic flows $\bold{U}^0$, the Euler correctors solve a linear hyperbolic system.  Then by a linear Prandtl equation, we can construct the weak boundary correctors $(u_p^i,v_p^i,\rho_p^i)$ to adjust the velocity to the non-slip boundary condition on the rigid walls.    Here we assume that the basic flows $\mathbf{U}^0$ are any shear flows satisfying (\ref{mu}) and generally speaking, if $\f{d^k\mu(x_2)}{dx_2^k}\neq0$ for $x_2=0,2,\ k=2,3$, then the 2nd-order compatibility condition of the hyperbolic system (\ref{e1}) will fail,  limiting the regularity of Euler correctors $(u_e^1,v_e^1,\rho_e^1)$ to $W^{2,\infty}(\Omega)$ globally.   On the other hand, by the classical theory of first-order linear hyperbolic system, the Euler correctors are piecewise smooth in domains separated by the characteristics which will be sufficient to ensure that $v_p^2\in W^{2,p}(\Omega)$. Finally, the boundary conditions of the  correctors $u_p^1,u_p^2,v_e^2,\rho_e^2$ on the inflow part of the boundary are properly chosen to satisfy the compatibility conditions  on the corner $(0,0)$, guaranteeing sufficient regularity for the approximate solutions: $u_s\in W^{3,p}(\Omega),\  v_s\in W^{2,p}(\Omega)$.

Then we  study the linearized system of (\ref{0.3})-(\ref{0.5}) around the approximate solution $(\rho_s,u_s,v_s)$, which is a hyperbolic-elliptic mixed system for the remainders $(\rho,u,v)$. The  control of  the total energy and total mass of the remainders are achieved by   exploiting  the supersonic property of the basic flows to perform a weighted estimate  which is the  key step in proving the existence of the solutions to the linearized system.  Based on the weighted estimate and  the classical method of  energy estimates, we can close the estimates of total energy and total mass, i.e. the estimate  (\ref{energyes}) holds and the existence of weak solutions to the linearized system follows immediately.

Finally, to prove the existence of solutions to the nonlinear system, we  need higher regularity of the remainders. Compared with the incompressible flows, the challenge in the proof of higher regularity for compressible flows lies in the perturbation of the density in the convection term and the mass equation.
In   the compressible case in 2-D, the nonlinear terms  are not fully controlled by $ H^2\times H^2\times H^1$ regularity.  What's more, considering the Dirichlet boundary condition of the velocity on $\partial\Omega$ and the Dirichlet boundary condition of the density on the inflow part of the boundary, it seems that we can hardly obtain that $(u,v,\rho)\in H^3\times H^3\times H^2$, and the $W^{2,p}\times W^{2,p}\times W^{1,p}$- estimate with $p>2$ seem to be a good choice. To prove the $W^{2,p}\times W^{2,p}\times W^{1,p}$ estimates of the remainders, our first important observation is that the density on the boundary of the domain can be well controlled. Then we consider the momentum equation as an elliptic system of the velocity. The main novelty in this step lies in the construction of a function $W=(W^1,W^2)$ defined in (\ref{lamme.3}) that satisfies an inhomogeneous elliptic system with homogeneous Dirichlet boundary condition on the boundary of the first quadrant. To construct $W$, we split it into two parts: $W_1=(W_{11}, W_{12})$ defined in (\ref{w1}) which satisfy an inhomogeneous elliptic system, and $W_2=(W_{21},W_{22})$ defined in (\ref{w2}),(\ref{w}) to adjust the boundary value of $W$ to the homogenous Dirichlet boundary value. In fact, by homogenizing the boundary value of $\rho$ to a new function $\hat \rho$, we can reduce the expression of $W_1$ to the convolution of $\nabla\hat\rho$ and the fundamental solution of Laplace operator. Then we can construct $W_2$ by use of the  the Green's function of Laplace operator in the first quadrant. By  Calderon-Zygmund theory and a detailed analysis we can prove that
$$\va\|W^1\|_{2,p;Q_{2R}(0)}+\va\|\partial_1W^2\|_{1,p;Q_{2R}(0)}\leq C\|\partial_1\hat\rho\|_{p;Q_{2R}(0)}.$$
Finally  by a careful bootstrap argument, the interpolation inequalities and the $L^p$ theory of elliptic systems in nonsmooth domains we can obtain the uniform-in-$\va$ $(W^{2,p})^2\times W^{1,p}$ estimates  of the remainders $(\rho, u,v)$.

The paper is organized as following: In section 2, we give the formal asymptotic expansion of $(u^\va,v^\va,\rho^\va)$ around the basic shear flows. First we  construct compressible Euler correctors $(u_e^i,v_e^i,\rho_e^i)$ with slip boundary condition on the rigid walls. Under the assumption that the the basic flow  is supersonic, the Euler correctors satisfy a linear hyperbolic system.  Then we will construct the weak boundary corrector $(u_p^i,v_p^i,\rho_p^i)$ to adjust the velocity to the non-slip boundary conditions on the rigid walls. In section 3, we will study the linearized system. First we prove the existence of solutions to a approximate system in Hilbert space. Then to deal with the nonlinear system, we also prove the $(W^{2,p})^2\times W^{1,p}$ estimates of solutions to the linear system. Finally, in section 4, we prove the existence of solutions to the nonlinear system (\ref{0.3})-(\ref{0.5}).

Now let us introduce the notations used throughout this paper.
\\[1mm]
{\sc Notation:} \ Let $G$ be an open set in $\mathbb{R}^N$.
We denote by $L^p(G)$ ($p\geq 1$) the Lebesgue spaces, by $W^{s,p}(G)$ ($p\geq 1$) the Sobolev spaces with $s$ being a real number,
by $H^k(G)$ ($k\in \mathbb{N}$) the Sobolev spaces $W^{k,p}(G)$ with $p=2$, and
by $C^k(G)$ (resp. $C^k(\overline{G})$) the space of $k$th-times continuously differentiable functions in $G$ (resp. $\overline{G}$).
 We use $|\cdot|_{k,p}$ to denote the standard norm in $W^{k,p}$ at the boundary $\partial\Omega$ and $|\cdot|$ for the norm in $L^2(G)$   throughout this paper.
 $\|\cdot\|_{k,p}$ stands for the standard norm in $W^{k,p}(G)$ and $\|\cdot\|$ for the norm in $L^2(G)$.
 We also use $\|\cdot\|_{L^\infty}$ to denote $\|\cdot \|_{L^\infty(\Omega)}=\text{ess~sup}_{\Omega}|\cdot |$.
 The symbol $\lesssim$ means that the left side is less than the right side multiplied by some constant.


We also define a smooth cut-off function $\chi(t)\in C^\infty([0,\infty))$ satisfying $|\chi|\leq1,\ |\chi|_{C^4([0,\infty))}\leq C$ and
\begin{equation}\label{cutoff}\chi(t)=\begin{cases}1,\ 0\leq t\leq \f34,\\0,\ t\geq1,\end{cases}\end{equation}
here $C>0$ is a finite constant.

\renewcommand{\theequation}{\thesection.\arabic{equation}}
\setcounter{equation}{0}

\section{Formal asymptotic expansion around shear flows}

In this section we will expand the solutions of the nonlinear system around the basic flows $\bold{U}^0 = (\mu(x_2), 0)$.  Here $\bold{U}^0$ satisfy the no-slip boundary condition on the moving boundaries and there is no mismatch between the basic flows and the moving boundaries. Thus, there would be no strong boundary layers around  $x_2=0$ and $x_2=2$. However, weak boundary layers will still arise due to viscous effects.
First we construct the  compressible Euler correctors
$(u_e^i,v_e^i,\rho_e^i)$ with slip boundary conditions  to balance the shear stress perturbations. Under the supersonic assumption of the basic flows, the Euler correctors solve a linear hyperbolic system.  Then  we  construct the weak boundary correctors $(u_p^i,v_p^i,\rho_p^i)$ to adjust the velocity to the non-slip boundary condition on the rigid walls.
 In what follows,
 the Eulerian profiles are functions of $(x_1,x_2)$, whereas the boundary layer profiles are functions of $(x_1,Y)$, where
\begin{align}
Y =
\left\{
\begin{aligned} \label{Y:defn}
&Y^+:= \frac{2-x_2}{\eps^{\frac 1 2}} \ \ \text{ if } 1 \le x_2 \le 2, \\
&Y^-:=\frac{x_2}{\eps^{\frac 1 2}}\ \  \text{ if } 0 \le x_2 \le 1.
\end{aligned}
\right.
\end{align}
Due to this, we break up the boundary layer profiles into two components, i.e.:
\begin{align}
u_p^i=
\left\{
\begin{aligned} \label{bp}
&u_p^{i,+}(x_1,Y^+) \ \ \text{ if } 1 \le x_2 \le 2, \\
&u_p^{i,-}(x_1,Y^-)\ \  \text{ if } 0 \le x_2 \le 1.
\end{aligned}
\right.
\end{align}
Then we expand the solutions in $\va$ as in (\ref{app.0}). In the following sections we will construct the Euler correctors and the weak boundary layer correctors separately.

\subsection{Euler correctors }

In this subsection we will focus on the construction of Euler correctors. The equations satisfied by the first Euler correctors are obtained by collecting the $\bigO(\eps)$  order Euler terms from (\ref{0.3})-(\ref{0.5}), and are  shown as following:
\begin{equation}
\begin{cases}
 \p_{x_1} u_e^1 + \p_{x_2} v_e^1 +\mu\p_{x_1}\rho_e^1= 0, \\
  \mu \p_{x_1} u_e^1 + \mu' v_e^1 + c^2\p_{x_1}  \rho_e^1 = \mu''(x_2) \\
  \mu \p_{x_1} v_e^1 + c^2\p_{x_2} \rho_e^1 = 0,
\end{cases} \label{e1}
\end{equation}
with the following boundary conditions:
$$ v_e^1|_{x_1 = 0} =u_e^1|_{x_1=0}=\rho_e^1|_{x_1=0}= 0,  v_e^1|_{x_2=0,2} = 0.$$
Similarly by collecting the $\bigO(\eps^{\f32})$  order Euler terms from (\ref{0.3})-(\ref{0.5}) we have
\begin{equation}\label{el.1}
\begin{cases}
 \p_{x_1} u_e^2 + \p_{x_2} v_e^2 +\mu\p_{x_1}\rho_e^2= 0, \\
  \mu \p_{x_1} u_e^2 + \mu' v_e^2 + c^2\p_{x_1}  \rho_e^2 = 0 \\
  \mu \p_{x_1} v_e^2 + c^2\p_{x_2} \rho_e^2 = 0,
\end{cases}
\end{equation}
with the following boundary conditions:
\begin{eqnarray*}&&u_e^2|_{x_1=0}=0,  v_e^2|_{x_2=0} = -v_p^1(x_1,0), v_e^2|_{x_2=2}=-v_p^1(x_1,2),
 v_e^2|_{x_1 = 0} =v^0(x_2),\\
 &&  \rho_e^2|_{x_1=0}= \rho^0(x_2).\end{eqnarray*}
Here the boundary conditions for $v_e^2$ on $x_2=0,2$ are chosen to adjust to the no-slip boundary condition, while the value of  $\rho_e^2,\ v_e^2$ on $x=0$ are properly chosen so that the first-order compatibility conditions on the corners are satisfied.
 More precisely we define
\begin{equation}
v^0(x_2)=\begin{cases}-v_p^1(0,0)\chi(\f {x_2}{b}),\ 0<x_2<1,\\-v_p^1(0,2)\chi(\f{2-x_2}b), \ 1\leq x_2<2,\end{cases}
\end{equation}
and
\begin{equation}
\rho^0(x_2)=\begin{cases}\f{\mu}{c^2}v_{px_1}^1(0,0)x_2\chi(\f {x_2}{b}),\ 0<x_2<1,\\-\f{\mu}{c^2}v_{px_1}^1(0,2)(2-x_2)\chi(\f{2-x_2}b), \ 1\leq x_2<2,\end{cases}
\end{equation}
here $b>0$ is a small constant.

If we denote by
\begin{align*}A=\left(
  \begin{array} {ccc}
    1 & 0 & \mu \\
   \mu & 0 & c^2 \\
   0 & \mu & 0\\
  \end{array}
\right),\ B=\left(
  \begin{array} {ccc}
    0 & 1 & 0 \\
   0 & 0 & 0 \\
   0 & 0 & c^2\\
  \end{array}
\right),\ D=\left(
  \begin{array} {ccc}
    0 & 0 & 0 \\
   0 & \mu' & 0 \\
   0 & 0 & 0\\
  \end{array}
\right),\ U^i=\left(\begin{array}{c}
u_e^i\\
v_e^i\\
\rho_e^i\\
 \end{array}\right),\end{align*}
 then system (\ref{e1}) and (\ref{el.1}) can be written as
 \begin{equation}
 AU_{x_1}^i+BU_{x_2}^i+DU^i=F_i,\ i=1,2,\ \ F_1=\left(\begin{array}{c}0\\ \mu'' \\ 0\\ \end{array}\right), \ F_2=\left(\begin{array}{c}0\\ 0 \\ 0\\ \end{array}\right).\label{s2}
 \end{equation}
 By solving det$(B-\la A)=0$, we find the eigenvalues of (\ref{e1}) and (\ref{el.1}) are:
 $$\la_1=0,\ \la_{2}=\f{c}{\sqrt{\mu^2-c^2}},\ \ \la_{3}=-\f{c}{\sqrt{\mu^2-c^2}}.$$
Obviously we have
$$-\f{c}{\sqrt{\mu^2-c^2}}<0<\f{c}{\sqrt{\mu^2-c^2}}.$$
 For $y_0\in[0,2]$, we define the i-th characteristics $y_i$ passing through $(0, y_0)$ by the ordinary differential equation
 \begin{align}
 &\f{dy_i(x_1;y_0)}{dx_1}=\la_i(y_i)\\
 &y_i(0;y_0)=y_0.
 \end{align}
For $\delta>0$, we denote by $\Omega_\delta=(0,\delta)\times(0,2)$ and
\begin{align}
&\Omega_1=\{(x_1,x_2)|0<x_1<\delta,\ 0<x_2<y_2(x_1;0)\}\\
&\Omega_{2}=\{(x_1,x_2)|0<x_1<\delta,\ y_2(x_1;0)<x_2<y_3(x_1,2)\}\\
&\Omega_{3}=\{(x_1,x_2)|0<x_1<\delta,\ y_3(x_1,2)<x_2<2\}.
\end{align}
First it is easy to check that for $i=1,2$, we have the following compatibility conditions on the corners $(0,0)$ and $(0,2)$:
$$\lim_{x_2\rightarrow0}v_e^i(0,x_2)=\lim_{x_2\rightarrow2}v_e^i(0,x_2)=\lim_{x_1\rightarrow0}v_e^i(x_1,0)=\lim_{x_1\rightarrow0}v_e^i(x_1,2)=0.$$
Moreover, from equation $(\ref{e1})_3$ and $(\ref{el.1})_3$  we can check that  the first-order compatibility conditions on the corners are satisfied, i.e.:
$$\lim_{x_2\rightarrow0}v_{ex_1}^i(0,x_2)=\lim_{x_1\rightarrow0}v_{ex_1}^i(x_1,0),\ \lim_{x_2\rightarrow2}v_{ex_1}^i(0,x_2)=\lim_{x_1\rightarrow0}v_{ex_1}^i(x_1,2).$$
By the theory of linear hyperbolic system (c.f.  Chapter 7 in \cite{L-Y}) we have:
\begin{theorem}\label{euler}
Assume  that $\mu(x_2)\in C^6([0,2])$ satisfying (\ref{mu}), $v_p^1$ is defined in (\ref{p.1}), then  there exists a constant $\delta>0$ such that for $i=1,2$, system (\ref{e1}) and (\ref{el.1}) has a unique solution $(u_e^i,v_e^i,\rho_e^i)\in C^{1,1}(\Omega_\delta)$   in $\Omega_\delta$ and the following estimates hold:
\begin{align}\label{asy.5}
\|\mau_e^i\|_{W^{2,\infty}(\Omega_\delta)}+\|\rho_e^i\|_{W^{2,\infty}(\Omega_\delta)}\leq C|\mu|_{C^6([0,2])},
\end{align}
Moreover, $(u_e^i,v_e^i,\rho_e^i)$ are piecewise smooth in $\Omega_j$($j=1,2,3$) with the following estimates:
\begin{equation}
\sum_{j=1}^3(\|\mau_e^1\|_{C^5(\bar{\Omega}_j)}+\|\rho_e^1\|_{C^5(\bar{\Omega}_j)})+\sum_{j=1}^3(\|\mau_e^2\|_{C^3(\bar{\Omega}_j)}+\|\rho_e^2\|_{C^3(\bar{\Omega}_j)})\leq C|\mu|_{C^6([0,2])}.\label{euler.1}
\end{equation}

\end{theorem}
\begin{remark}
From system (\ref{e1}), it is easy to check that
$$\partial_{x_1}\rho_e^1|_{x_1=0}=\f{\mu''}{c^2-\mu^2},\ \partial_{x_1x_2}\rho_e^1|_{x_1=0}=\f{\mu'''}{c^2-\mu^2}+\f{2\mu\mu'\mu''}{(c^2-\mu^2)^2},$$
while from $(\ref{e1})_3$ we have
$$\lim_{x_1\rightarrow0}\partial_{x_1x_2}\rho_e^1(x_1,0)=\lim_{x_1\rightarrow0}\partial_{x_1x_2}\rho_e^1(x_1,2)=0.$$
Generally speaking, if $\mu''(0)\neq0$ or $\mu''(2)\neq0$, then we can not obtain the compatibility conditions for $\partial_{x_1x_2}\rho_e^1$ on corners $(0,0)$ or $(0,2)$. So for general shear flow $\mathbf{U}_0$ satisfying (\ref{mu}), we can only obtain the global $C^{1,1}$ regularity in Theorem \ref{euler}.
\end{remark}

\subsection{Weak boundary layer correctors }

In this subsection we will construct the  weak boundary layer correctors to adjust the velocity to the non-slip boundary conditions on the rigid boundaries.
First the leading $\bigO(\eps^{\f12})$ order boundary layer terms from (\ref{0.3})-(\ref{0.5}) is:
$$\p_{Y} \rho_p^1 = 0,$$
and the leading $\bigO(\eps)$ order boundary layer terms  are:
$$u_{px_1}^1+v_{pY}^1=0\ \ \ \text{and}\ \ \  \mu \p_{x_1} u_p^1 - \p_{YY}u_p^1 +c^2\partial_{x_1}\rho_p^1\triangleq R^{u,1},$$
As we will assume that the boundary layer correctors $(u_p^i,v_p^i,\rho_p^i)$  decrease rapidly to $0$ when $Y$ tends to $\infty$, we have obviously  that $\rho_p^1\equiv0$. To construct $(u_p^1,v_p^1)$ , we first consider the following initial-boundary value problem of parabolic equation with constant coefficients:
\begin{eqnarray} \label{pr.BVP.1}
\begin{cases}
A^- \p_{x_1} u^{1,0}_p - \p_{YY}u^{1,0}_p = 0, \\
u^{1,0}_p|_{x_1 = 0} =-(\sum_{k=1}^4\f1{(2k)!}(A^-)^k\partial_{x_1}^ku_e^1(0,0)Y^{2k})\chi(Y) ,\\
 u^{1,0}_p|_{Y = 0} = -u_e^1|_{x_2 = 0}, u^{1,0}_p|_{Y \rightarrow +\infty} = 0\\
 v^{1,0 }_p = \int_{Y}^{\infty}\p_{x_1} u^{1,0 }_p.
\end{cases}
\end{eqnarray}
here $A^-=\mu(0)$, $\chi(Y)$ is defined in (\ref{cutoff}). It is easy to check that the fourth-order compatibility conditions on the corner $(0,0)$ are satisfied. If we denote by $\Omega_p=(0,L)\times (0,+\infty)$, then  by Lemma \ref{ba}, there exists a unique solution $(u_p^{1,0},v_p^{1,0})\in W^{5,p}(\Omega_p)\times W^{4,p}(\Omega_p)$ satisfying system (\ref{pr.BVP.1}) and for $w(Y)$ defined in (\ref{A1}), $ 0 \le 2k+ l \le 10, m, j\in N$, we have the following estimate:
\begin{eqnarray*}
&&\|(1+Y)^m w(Y)\nabla ^j\mau_p^{1,0} \|_{L^\infty}+\|\partial_{x_1}^k\partial_Y^{l}\mau_p^{1,0}\|_{p}\lesssim|u_e^1(\cdot,0)|_{c^5([0,L])}.
\end{eqnarray*}

Then we cut-off $\mau_p^{1,0}=(u_p^{1,0},v_p^{1,0})$ to obtain the first  boundary layer correctors $\mau_p^{1,-}=(u_p^{1,-},v_p^{1,-})$ near $x_2=0$:
\begin{align} \label{cut.off.2}
u^{1,-}_p = \chi(\frac{\sqrt{\eps}Y^-}{a_0}) u^{1,0}_p - \frac{\sqrt{\eps}}{a_0} \chi'(\frac{\sqrt{\eps}Y^-}{a_0}) \int_0^{x_1} v_p^{1,0}, \hspace{2 mm} v^{1,-}_p := \chi(\frac{\sqrt{\eps}Y^-}{a_0}) v^{1,0}_p,
\end{align}
where $a_0>0$ is a constant small enough.
After cutting off (\ref{cut.off.2}), we have the  contribution with $O(\va^{\f12})$ order to the remainders
\begin{eqnarray*}
\mathcal{C}_{cut}^{1,-}=&&\f{A^-}{a_0}\eps^{\f12}\chi'v_p^{1,0}-\f2{a_0}\eps^{\f12}\chi'\partial_{Y^-}u_p^{1,0}
-\f2{a_0^2}\eps\chi''u_p^{1,0}-\f{\eps^{\f32}}{a_0^3}\chi'''\int_0^{x_1}v_p^{1,0}(s,Y^-)ds,
\end{eqnarray*}
i.e. we have
\begin{eqnarray}
\begin{cases}
A^- \p_{x_1} u^{1,-}_p - \p_{YY}u^{1,-}_p =\mathcal{C}_{cut}^{1,-} , \\
u^{1,-}_p|_{x_1 = 0} =-(\sum_{k=1}^4\f1{(2k)!}(A^-)^k\partial_{x_1}^ku_e^1(0,0)Y^{2k})\chi(Y),\\
 u^{1,-}_p|_{Y = 0} = -u_e^1|_{x_2 =0}, v^{1,-}_p|_{Y=0}=0,\ u^{1,-}_p|_{Y \rightarrow +\infty} = 0\\
 u^{1,-}_{px_1}+v^{1,- }_{pY} = 0.
\end{cases}
\end{eqnarray}

Besides, due to the approximation of $\mu(x_2)$ by $\mu(0)$ in the support of the cut-off function $\chi(\frac{\eps^{\frac 12} Y^-}{a_0})$, we have another contribution with $O(\va^{\f12})$ order to the error defined by
\begin{align}
\mathcal{C}_{app}^{1,-} := &  (\mu(x_2)- \mu(0) )[\chi(\frac{\va^{\frac 12} Y^-}{a_0}) \p_{x_1} u^{1,0}_p-\f1{a_0}\eps^{\f12}\chi'v_p^{1,0}],
\end{align}
i.e. we have
$$\mu \p_{x_1} u^{1,-}_p - \p_{YY}u^{1,-}_p =\mathcal{C}_{cut}^{-}+\mathcal{C}_{app}^{1,-}.$$
The constructions of  $\mau_p^{1,+},\mathcal{C}_{cut}^{1,+},\mathcal{C}_{approx}^{1,+}$ near the boundary $x_2=2$ are exactly the same as above. If we denote by
$$u_p^1=u_p^{1,+}+u_p^{1,-}, \ v_p^1=v_p^{1,+}+v_p^{1,-},\ \rho_p^1\equiv0,\  \mathcal{C}_{cut}^1=\mathcal{C}_{cut}^{1,+}+\mathcal{C}_{cut}^{1,-},\ \mathcal{C}_{app}^1=\mathcal{C}_{app}^{1,-}+\mathcal{C}_{app}^{1,+}.$$
and
\begin{equation}u^{1,0}(Y)=\begin{cases}
-(\sum_{k=1}^4\f1{(2k)!}(\mu(0))^k\partial_{x_1}^ku_e^1(0,0)Y^{2k})\chi(Y),\ \ 0\leq x_2\leq 1,\\
-(\sum_{k=1}^4\f1{(2k)!}(\mu(2))^k\partial_{x_1}^ku_e^1(0,2)Y^{2k})\chi(Y),\ \  1\leq x_2\leq 2,
\end{cases}\label{pb}\end{equation}
then the first boundary layer correctors $(u_p^1,v_p^1)$ satisfy the following system:
\begin{eqnarray}
\begin{cases}
\mu \p_{x_1} u_p^1 - \p_{YY}u_p^1 =\mathcal{C}_{cut}^1+\mathcal{C}_{app}^1 , \\
u_p^1|_{x_1 = 0} = u^{1,0}(Y),\   u_p^1|_{x_2 = 0} = -u_e^1(x_1,0), u_p^1|_{x_2 = 2} = -u_e^1(x_1,2),\\
u_{px_1}^1+v_{pY}^1=0
\end{cases}\label{p.1}
\end{eqnarray}
To construct  the second boundary layer correctors, we will consider the following parabolic problem for $u_p^{2,0}$:
\begin{eqnarray} \label{pr.BVP.2}
\begin{cases}
A^- \p_{x_1} u^{2,0}_p - \p_{YY}u^{2,0}_p = 0, \\
u^{2,0}_p|_{x_1 = 0} =-[\f12A^-\partial_{x_1}u_e^2(0,0)Y^2+\f1{4!}(A^-)^2\partial_{x_1x_1}u_e^2(0,0)Y^4]\chi(Y) ,\\
 u^{2,0}_p|_{Y = 0} = -u_e^2|_{y = 0}, u^{2,0}_p|_{Y \rightarrow +\infty} = 0,
\end{cases}
\end{eqnarray}
while $ v^{2,0 }_p$ is defined as
$$v^{2,0 }_p(x_1,Y) = -\int_0^{Y}\p_{x_1} u^{2,0 }_p(x_1,s)ds.$$
Then we cut-off $\mau_p^{2,0}=(u_p^{2,0},v_p^{2,0})$ to obtain the second boundary layer corrector $\mau_p^{2,-}=(u_p^{2,-},v_p^{2,-})$ near $x_2=0$:
\begin{align} \label{cut.off.3}
u^{2,-}_p = \chi(\frac{\sqrt{\eps}Y^-}{a_0}) u^{2,0}_p - \frac{\sqrt{\eps}}{a_0} \chi'(\frac{\sqrt{\eps}Y^-}{a_0}) \int_0^{x_1} v_p^{2,0}, \hspace{2 mm} v^{2,-}_p := \chi(\frac{\sqrt{\eps}Y^-}{a_0}) v^{2,0}_p.
\end{align}
Similarly we can construct $\mau_p^2=(u_p^2,v_p^2)$ and $\mathcal{C}_{cut}^2, \mathcal{C}_{app}^2$ satisfying
\begin{eqnarray}\begin{cases}
\mu \p_{x_1} u_p^2 - \p_{YY}u_p^2 =\mathcal{C}_{cut}^2+\mathcal{C}_{app}^2 , \\
u_p^2|_{x_1 = 0} = u^{2,0}(Y),\   u_p^2|_{x_2 = 0} = -u_e^2(x_1,0), u_p^2|_{x_2 = 2} = -u_e^2(x_1,2),\\
 v_p^2|_{x_2=0}=v_p^2|_{x_2=2}=0,\\
u_{px_1}^2+v_{pY}^2=0
\end{cases}\label{p.2}
\end{eqnarray}
with
\begin{equation}u^{2,0}(Y)=\begin{cases}
-[\f12\mu(0)\partial_{x_1}u_e^2(0,0)Y^2+\f1{4!}\mu(0)^2\partial_{x_1x_1}u_e^2(0,0)Y^4]\chi(Y),\ \ 0\leq x_2\leq 1,\\
-[\f12\mu(2)\partial_{x_1}u_e^2(0,2)Y^2+\f1{4!}\mu(2)^2\partial_{x_1x_1}u_e^2(0,2)Y^4]\chi(Y),\ \  1\leq x_2\leq 2.
\end{cases}\end{equation}
Finally, the main result of this subsection reads as:

\begin{theorem}\label{boundarylayer}
For given $(u_e^i,v_e^i,\rho_e^i)$, $i=1,2$, defined in Theorem \ref{euler},  there exists a solution $(u_p^1,v_p^1)\in W^{5,p}(\Omega)\times W^{4,p}(\Omega)$ satisfying system (\ref{p.1}) and  a solution $(u_p^2,v_p^2)\in W^{3,p}(\Omega)\times W^{2,p}(\Omega)$ satisfying system (\ref{p.2}). Besides, for $0\leq l\leq 2$, $m,k,j,j_1,j_2,k_1,k_2\geq0$, $0\leq2k_1+j_1\leq10$, $0\leq2k_2+j_2\leq6$, we have the following estimates:
\begin{eqnarray}
&&\|(1+Y)^m w(Y)\nabla ^j\mau_p^1 \|_{L^\infty}+\|\partial_{x_1}^{k_1}\partial_Y^{j_1}\mau_p^1\|_{L^p(\Omega_p)}+\|(1+Y)^m w(Y)\nabla ^jv_{pY}^2 \|_{L^\infty}\nonumber\\
&&+\|(1+Y)^m w(Y)\nabla ^ju_p^2 \|_{L^\infty}+\|\partial_{x_1}^{k_2}\partial_Y^{j_2}u_p^2\|_{L^p(\Omega_p)}+\|\partial_{x_1}^lv_p^2\|_{L^\infty(\Omega_p)}\nonumber\\
\leq&&C|\mu|_{c^6([0,2])}.\label{p.4}
\end{eqnarray}
 Moreover, we have
 \begin{eqnarray}
&&\sum_{i=1}^2\|\mathcal{C}_{cut}^i\|_{L^\infty(\Omega)}+\sum_{i=1}^2\|\mathcal{C}_{app}^i\|_{L^\infty(\Omega)}
\leq C\va^{\f12}|\mu|_{c^6([0,2])},\label{r.2}
\end{eqnarray}
and
\begin{eqnarray}
&&\|\mathcal{C}_{cut}^1\|_{L^p(\Omega)}+\|\mathcal{C}_{app}^1\|_{L^p(\Omega)}
\leq C\va^{\f12+\f1{2p}}|\mu|_{c^6([0,2])}\label{r.1}
\end{eqnarray}
where the constant $C$ is independent of $\va$.
\end{theorem}
\begin{proof}
First of all, the existence follows directly from the process above.
By Lemma \ref{ba}, to prove (\ref{p.4}), we only need to prove the estimates of $\|\partial_{x_1}^lv_p^2\|_{L^\infty(\Omega_p)}$. In fact, for $0\leq k  \leq 2$,
\begin{eqnarray*}
&&\|\partial_{x_1}^kv_p^2\|_{L^\infty(\Omega_p)}\leq \sup_{0\leq x_1\leq L}\int_0^{\infty}|\partial_{x_1}^ku_{px_1}^2(x_1,s)|dY\\
\leq&&\|(1+Y)^2 w(Y)\nabla ^3u_p^2 \|_{L^\infty(\Omega_p)}\int_0^{\infty}(1+Y)^{-2}dY+\|u_p^2\|_{3,p;\Omega_p} \\
\lesssim &&|\mu|_{c^6([0,2])}.
\end{eqnarray*}
Next estimate (\ref{r.2}) follows immediately from (\ref{p.4}) and the definition of $\mathcal{C}_{cut}^i,\mathcal{C}_{app}^i$.
Finally, by (\ref{p.4}) and direct computation we have
\begin{eqnarray*}
&&\|\mathcal{C}_{cut}^1\|_{L^p(\Omega)}^p+\|\mathcal{C}_{app}^1\|_{L^p(\Omega)}^p\nonumber\\
\lesssim &&\va^2|\mau_p^1|_{L^\infty}^2+\va^{\f p2}\int_\Omega(|v_p^1|+|u_{pY}^1|+Y|u_{px_1}^1|)^pdx_1dx_2\\
\lesssim &&\va^2|\mau_p^1|_{L^\infty}^2+\va^{\f{p+1}2}\int_{\Omega_p}(|v_p^{1,\pm}|+|u_{pY}^{1,\pm}|+Y|u_{px_1}^{1,\pm}|)^pdx_1dY\\
\lesssim &&\va^{\f{p+1}2}|\mu|_{c^6([0,2])},
\end{eqnarray*}
and we have finished the proof.
\end{proof}

\subsection{The approximate solution $(u_s,v_s,\rho_s)$}

Recall that
\begin{eqnarray*}
&&u_s=\mu+\va u_e^1+\va u_p^1+\va^{\f32}u_e^2+\va^{\f32}u_p^2,\\
&&v_s=\va v_e^1+\va^{\f32}v_p^1+\va^{\f32}v_e^2+\va^2v_p^2,\\
&&\rho_s=\rho^*+\va \rho_e^1+\va^{\f32}\rho_e^2,
\end{eqnarray*}
we have by (\ref{e1}),(\ref{el.1}), (\ref{p.1}), (\ref{p.2}) and direct computation that
\begin{eqnarray}
&&\rho_s\dv\mau_s+ \mau_s\cdot\nabla\rho_s=g_{0s}\label{g1}\\
&&\rho_s\mau_s\cdot\nabla\mau_s-\va\Delta\mau-\la\va\nabla\dv\mau_s+c^2\nabla \rho_s=\mathbf{g}_s,\label{g2}
\end{eqnarray}
here $\mathbf{g}_s=(g_{1s},g_{2s})$ with
\begin{eqnarray*}
g_{0s}=&&\va^2(\rho_e^1+\va^{\f12}\rho_e^2)(\dv\mau_e^1+\va^{\f12}\dv\mau_e^2)+\va^2(\rho_e^1+\va^{\f12}\rho_e^2)_{x_1}(u_e^1+u_p^1+\va^{\f12}u_e^2+\va^{\f12}u_p^2)\\
&&+\va^2(\rho_e^1+\rho_e^2)_{x_2}(v_e^1+\va^{\f12}v_p^1+\va^{\f12}v_e^2+\va v_p^2)\\
g_{1s}=&&\va^{\f32}[a\mu'v_p^1+u_{pY}^1v_e^1]+\va\mathcal{C}_{cut}^1+\va\mathcal{C}_{app}^1+\va^{\f32}\mathcal{C}_{cut}^2+\va^{\f32}\mathcal{C}_{app}^2\\
&&+\va^2 u_s(\rho_e^1+\va^{\f12}\rho_e^2)(u_e^1+u_p^1+\va^{\f12}u_e^2+\va^{\f12}u_p^2)_{x_1}\\
&&+a\va^2(u_e^1+u_p^1+\va^{\f12}u_e^2+\va^{\f12}u_p^2)(u_e^1+u_p^1+\va^{\f12}u_e^2+\va^{\f12}u_p^2)_{x_1}\\
&&+\va^2[\rho_s\mu'v_p^2+(v_e^1+\va^{\f12}v_p^1+\va^{\f12}v_e^2+\va v_p^2)(u_e^1+u_p^1+\va^{\f12}u_e^2+\va^{\f12}u_p^2)_{x_2}]\\
&&-\va^2(\rho_e^1+\va^{\f12}\rho_e^2)(v_e^1+\va^{\f12}v_p^1+\va^{\f12}v_e^2+\va v_p^2)u_{sx_2}+\va^2(1+\la)( u_{ex_1x_1}^1+\va^{\f12}u_{ex_1x_1}^2)\\
&&+\va^2 (u_{ex_2x_2}^1+\va u_{ex_2x_2}^2+u_{px_1x_1}^1+\va^{\f12}u_{px_1x_1}^2)\\
g_{2s}&&=\va^{\f32}a\mu (v_{px_1}^1+\va^{\f12}v_{px_1}^2)+\va^{\f32}(1+\la)v_{pYY}\\
&&+\va^2u_s(\rho_e^1+\va^{\f12}\rho_e^2)(v_e^1+\va^{\f12}v_p^1+\va^{\f12}v_e^2+\va v_p^2)_{x_1}\\
&&+\va^2a(u_e^1+u_p^1+\va^{\f12}u_e^2+\va^{\f12}u_p^2)(v_e^1+\va^{\f12}v_p^1+\va^{\f12}v_e^2+\va v_p^2)_{x_1}\\
&&+\va^2\rho_s(v_e^1+\va^{\f12}v_p^1+\va^{\f12}v_e^2+\va v_p^2)(v_e^1+\va^{\f12}v_p^1+\va^{\f12}v_e^2+\va v_p^2)_{x_2}+\va^2(1+\la)v_{pYY}^2\\
&&+\va^{\f52}(v_{px_1x_1}^1+\va^{\f12}v_{px_1x_1}^2)+\va^2(\Delta (v_e^1+\va^{\f12}v_e^2)+\la \partial_{x_2}\dv(\mau_e^1+\va^{\f12}\mau_e^2)
\end{eqnarray*}

The main result of this section reads as:
\begin{theorem}\label{thm2.2}
Assume that $\mu$  is a smooth function satisfying conditions in Theorem  \ref{main}, then we can construct a triple  $(u_s,v_s, \rho_s)\in W^{2,p}(\Omega)\times W^{2,p}(\Omega)\times W^{2,\infty}(\Omega)$ with formula defined in (\ref{app.0}) and:
\begin{equation}\dv\mau_p^i=0\ \text{in}\ \Omega,\ u_s(0)=\mu(0),\ u_s(2)=\mu(2), v_s(0)=v_s(2)=0.\label{2.2.1}\end{equation}
Moreover, the Euler correctors $(u_e^i,v_e^i,\rho_e^i)$ satisfy estimates (\ref{asy.5}), (\ref{euler.1}) while the weak boundary layer correctors $(u_p^i,v_p^i)$ satisfy estimate (\ref{p.4})-(\ref{r.1}) and $\rho_p^i\equiv0$. Finally we have the following estimates:
\begin{eqnarray}
&&\|g_{0s}\|_{L^2}+\va^{\f12}\|g_{0s}\|_{W^{1,p}}\lesssim \va^2,\\
&&\|\mathbf{g}\|_{L^p}\lesssim \va^{\f32+\f1{2p}}.
\end{eqnarray}

\end{theorem}

The proof of Theorem \ref{thm2.2} follows immediately from Theorem \ref{euler} and Theorem \ref{boundarylayer}.

\bigskip

\renewcommand{\theequation}{\thesection.\arabic{equation}}
\setcounter{equation}{0}

\section{The Linearized System}

To prove Theorem \ref{main}, we will first study the linearized system of (\ref{0.3})-(\ref{0.5}) around the approximate solutions $(u_s,v_s,\rho_s)$ constructed above.
Putting the expansion (\ref{app.0}) into (\ref{0.3})-(\ref{0.5}), we find that the remainder solutions $(\rho,u,v)$ satisfying the following system:
\begin{eqnarray}
\text{div}\mathbf{u}+  u^\va\rho_{x_1}+   v^\va \rho_{x_2}&=&g_0(\rho,u,v)\qquad\text{in}~\Omega,\label{1.1}\\
u_s u_{x_1}+u_{sx_2}v-\va\Delta u
-\va\partial_{x_1}\mathrm{div}\mathbf{u}+c^2\rho_{x_1}&=&g_1(\rho,u,v)\qquad\text{in}~\Omega,\label{1.2}\\
u_s v_{x_1}-\va\Delta v
-\va\partial_{x_2}\mathrm{div}\mathbf{u}+c^2\rho_{x_2}&=&g_2(\rho,u,v)\qquad\text{in}~\Omega,\label{1.3}
\end{eqnarray}
here
\begin{eqnarray*}
&&g_0(\rho,u,v)=g_{0s}+g_{0r}(\rho,u,v),\\
&& g_1(\rho,u,v)=g_{1s}+g_{1r}(\rho,u,v),\\
&& g_2(\rho,u,v)=g_{2s}+g_{2r}(\rho,u,v)\end{eqnarray*}
and $g_{0s},\mathbf{g}_s$ are defined in (\ref{g1}), (\ref{g2}) while
\begin{eqnarray*}
g_{0r}(\rho,u,v)&&=-\va\rho\dv\mathbf{U}_e-\va\mau\cdot\nabla P_e-\va P_e\dv\mau-\rho\dv\mau\\
g_{1r}(\rho,u,v)&&
=-\rho_s(uu_{x_1}+vu_{x_2})-\va[(U_e+U_p)\rho u_{x_1}+(U_e+U_p)_{x_1}\rho u]\\
&&-\va\rho_s[u(U_{ex_1}+U_{px_1})+(V_e+\va^{\f12}V_p)u_{x_2}]-\va P_eu_su_{x_1}-\va\rho\mu(U_e+U_p)_{x_1}\\
&&-\va\mu'\rho(V_e+\va^{\f12}V_p)-\mu\rho u_{x_1}-\va^2\rho(U_e+U_p)(U_e+U_p)_{x_1}-\rho uu_{x_1}\\
&&+\rho_svu_{x_2}+\va\rho u_{x_2}(V_e+\va^{\f12}V_p)-\va P_e\mu'v-\mu'\rho v-\va(\rho_s+\rho)v U_{ex_2}\\
&&-\va^{\f12}(\rho_s+\rho)v U_{pY}-\rho vu_{x_2}+[c^2-p'(\rho^\va)]\rho_{x_1}\\
g_{2r}(\rho,u,v)&&
=-\va[P_eu_sv_{x_1}+(\rho_su+\rho u_s+\rho u)(V_e+\va^{\f12}V_p)_{x_1}]-(\rho_su+\rho u_s)v_{x_1}\\
&&-\va[\rho_s(V_e+\va^{\f12}V_p)v_{x_2}+\rho_s v(V_{ex_2}+V_{pY})+\rho(V_e+\va^{\f12}V_p)v_{x_2}+\rho vV_{ex_2}]\\
&&-\va\rho vV_{pY}-\va^2(V_e+\va^{\f12}V_p)(V_{ex_2}+V_{pY})-\rho_svv_{x_2}-\rho vv_{x_2}\\
&&-\rho uv_{x_1}+[c^2-p'(\rho^\va)]\rho_{x_2},
\end{eqnarray*}
here $\mathbf{U}_e=(U_e,V_e)$, $\mathbf{U}_p=(U_p,V_p)$ and
\begin{eqnarray*}U_e=u_e^1+\va^{\f12}u_e^2,\ V_e=v_e^1+\va^{\f12}v_e^2,\\
U_p=u_p^1+\va^{\f12}u_p^2,\ V_p=v_p^1+\va^{\f12}v_p^2, P_e=(\rho_e^1+\va^{\f12}\rho_e^2)
\end{eqnarray*}
To remove the inhomogeneity from the boundary conditions, we define
$$\bar u=u-b_1,\ \bar v=v-b_2,\ \bar\rho=\rho-(\rho_0(x_2)-\rho_s(0,x_2)),$$
here $\mathbf{b}=(b_1,b_2)$ and
\begin{eqnarray*}
&&b_1(x_1,x_2)=(1-\frac {x_1}L)[u_{in}-u_s(0,x_2)]+\f {x_1}L[u_{out}-u_s(L,x_2)]\\
&&b_2(x_1,x_2)=(1-\frac {x_1}L)[v_{in}-v_s(0,x_2)]+\f {x_1}L[v_{out}-v_s(L,x_2)].
\end{eqnarray*}
Then $(\bar u,\bar v,\bar\rho)$ satisfy the following system
\begin{eqnarray}
\text{div}\bar \mau+   u^\va\bar \rho_{x_1}+    v^\va\bar  \rho_{x_2}&=&\bar g_0+g_0(\rho,u,v)\qquad\text{in}~\Omega,\label{linear.4}\\
u_s\bar  u_{x_1}+u_{sx_2}\bar v-\va\Delta\bar  u
-\va\partial_{x_1}\mathrm{div}\bar \mau+c^2\bar \rho_{x_1}&=&\bar g_1+g_1(\rho,u,v)\qquad\text{in}~\Omega,\\
u_s\bar  v_{x_1}-\va\Delta \bar v
-\va\partial_{x_2}\mathrm{div}\bar \mau+c^2\bar \rho_{x_2}&=&\bar g_2+g_2(\rho,u,v)\qquad\text{in}~\Omega,\label{linear.5}
\end{eqnarray}
with homogeneous boundary condition
\begin{eqnarray}
\bar\rho|_{x_1=0}=\bar u|_{\partial\Omega}=\bar v|_{\partial\Omega}=0,\label{boundary}
\end{eqnarray}
here
\begin{eqnarray*}
&&\bar g_{0s}=-b_{1x_1}-b_{2x_2}-v_s(\rho_{0}'-\rho_{sx_2}(0,x_2))\\
&&\bar g_0=\bar g_{0s}-v(\rho_{0}'-\rho_{sx_2}(0,x_2)),\\
&&\bar g_1=-u_sb_{1x_1}-u_{sx_2}b_2+\va\Delta b_1+\va\partial_{x_1}(b_{1x_1}+b_{2x_2}),\\
&&\bar g_2=-u_sb_{2x_1}+\va\Delta b_2+\va\partial_{x_2}(b_{1x_1}+a_{2x_2})-c^2(\rho_{0}'-\rho_{sx_2}(0,x_2)).
\end{eqnarray*}
For simplicity, we will omit the superscript in the following.
To prove the existence of solution to the nonlinear system (\ref{linear.4})-(\ref{boundary}), we will first construct a sequence from the following linear system that will converge
to a solution of the nonlinear system:
\begin{eqnarray}\label{n0}\begin{cases}
\text{div}\mathbf{u}^{n+1}+(\mau_s+\mathbf{b}+\mathbf{u}^{n})\cdot\nabla\rho^{n+1}
=g_0(\mathbf{u}^{n},\rho^{n})+\bar g_0(v^n)\ \text{in}~\Omega,\\
u_su^{n+1}_{x_1}+u_{sx_2}v^{n+1}-\va\Delta
u^{n+1}
-\va\lambda\partial_{x_1}\mathrm{div}\mathbf{u}^{n+1}+c^2\partial_{x_1}\rho^{n+1}\\
\qquad\qquad\qquad\qquad\qquad\qquad\qquad\qquad\qquad\quad=g_1(\mathbf{u}^{n},\rho^{n})+\bar g_1\ \text{in}~\Omega,\\
u_sv^{n+1}_{x_1}-\va\Delta
v^{n+1}
-\va\lambda\partial_{x_2}\mathrm{div}\mathbf{u}^{n+1}+c^2\partial_{x_2}\rho^{n+1}=g_2(\mathbf{u}^{n},\rho^{n})+\bar g_2\ \text{in}~\Omega,\\
\rho^{n+1}=0 ~ \text{on}~\Gamma_{\text{in}},\\
\mathbf{u}^{n+1}=0\ \text{on}~\partial\Omega,
\end{cases}\end{eqnarray}
where $g_0, g_1, g_2, \bar g_0, \bar  g_1, \bar g_2$ are defined  above, $c^2=a\gamma (\rho^*)^{\gamma-1}$ and $(\mathbf{u}^0,\rho^0)=(0,0,0)$.

In the following  we will first  consider  the linear system:
\begin{eqnarray}
\text{div}\mathbf{u}+ u^\va\rho_{x_1}+ v^\va \rho_{x_2}&=&\hat g_0\qquad\text{in}~\Omega,\label{linear1}\\
u_su_{x_1}+u_{sx_2}v-\va\Delta u
-\va\lambda\partial_{x_1}\mathrm{div}\mathbf{u}+c^2\rho_{x_1}&=&\hat g_1\qquad\text{in}~\Omega,\label{linear2}\\
u_sv_{x_1}-\va\Delta v
-\va\lambda\partial_{x_2}\mathrm{div}\mathbf{u}+c^2\rho_{x_2}&=&\hat g_2\qquad\text{in}~\Omega,\label{linear3}
\end{eqnarray}
with boundary condition (\ref{boundary}). Here $\mau^\va=(u^\va,v^\va)\in W^{2,p}(\Omega)$, $\hat g_0\in W^{1,p}(\Omega), \hat {\mathbf{g}}=(\hat g_1,\hat g_2)\in L^p(\Omega)$ are given functions satisfying $\|\mau^\va-\mau_s\|_{2,p;\Omega}\ll1$.

Before the proof we  straighten the stream line by introducing the following change of variables
$\Pi:\Omega\rightarrow\hat\Omega$
\begin{eqnarray}\begin{cases}x_1=\bar x_1,\\[2mm]
 \displaystyle{ x_2=\bar x_2+\int_0^{\bar x_1}\f{v^\va}{u^\va}(s,x_2(s,\bar x_2))ds.}
\end{cases}\label{2.79}\end{eqnarray}
Recalling that $v^\va=0$ on $\partial\Omega$ and $\|\mau^\va-\mau_s\|_{2,p;\Omega}$ is small enough, it is easily to check
that $\hat \Omega=[0,1]\times[0,1]$ and the mapping $\Pi$ is a diffeomorphism. A direct computation shows that
$$\|\f{\partial(\bar x_1,\bar x_2)}{\partial(x_1,x_2)}-I\|_{1,p;\Omega}\leq
C\|v^\va\|_{2,p;\Omega},$$
where $I$ is the unite matrix.
If we denote $\hat{\Gamma}_{\text{in}}=\Pi^{-1}(\Gamma_{\text{in}})$, $\hat{\Gamma}_{\text{out}}=\Pi^{-1}(\Gamma_{\text{out}})$, $\hat{\Gamma}_0=\Pi^{-1}(\Gamma_0)$, then the system (\ref{linear1})--(\ref{linear3}) can be rewritten as
(for convenience, we omit the superscript in the new coordinates)
\begin{eqnarray}\label{linear.0}\begin{cases}
\text{div}\mathbf{u}+u_s\partial_{x_1}\rho=f_0\qquad\text{in}~\Omega,\\
u_s\partial_{x_1}u+u_{sx_2}v-\va\Delta u-\va\lambda\partial_{x_1}\mathrm{div}\mathbf{u}+c^2\partial_{x_1}\rho=f_1\qquad\text{in}~\Omega,\\
u_s\partial_{x_1}v-\va\Delta v
-\va\lambda\partial_{x_2}\mathrm{div}\mathbf{u}+c^2\partial_{x_2}\rho=f_2\qquad\text{in}~\Omega,\\
\rho=0 ~\quad\text{on}~\Gamma_{\text{in}},\\
\mathbf{u}=0\qquad\text{on}~\partial\Omega,
\end{cases}\end{eqnarray}
where
\begin{eqnarray}
f_0&&=g_0+(1-\f{u_s}{u^\va})\dv\mau+\f{u_s}{u^\va}[g_0+u_{\bar x_1}(1-\f{\partial\bar x_1}{\partial x_1})+v_{\bar x_2}(1-\f{\partial\bar x_2}{\partial x_2})
+u_{\bar x_2}\f{\partial\bar x_2}{\partial x_1}]\nonumber\\
f_1&&=-u_s\partial_{\bar x_2}u\f{\partial\bar
x_2}{\partial x_1}+\va(1+\lambda)[2\partial_{\bar x_1\bar x_2}u\f{\partial\bar
x_2}{\partial x_1}+\partial_{\bar x_2}u\f{\partial^2\bar{x}_2}{\partial x_1^2}+\partial_{\bar
x_2}^2u(\f{\partial\bar x_2}{\partial x_2})^2]\nonumber\\
&& +\va\partial_{\bar
x_2}u\f{\partial^2\bar x_2}{\partial x_2^2}+\va[(\f{\partial\bar x_2}{\partial
x_2})^2-1]\partial_{\bar
x_2}^2u+\va\la\partial_{\bar x_1\bar x_2}v(\f{\partial\bar
x_2}{\partial x_2}-1)\nonumber\\
&&+\va\la(v_{\bar x_2\bar x_2}\f{\partial\bar x_2}{\partial x_2}\f{\partial\bar x_2}{\partial x_1}+v_{\bar x_2}\f{\partial^2\bar x_2}{\partial x_1\partial x_2})
-c^2\rho_{\bar x_2}\f{\partial\bar x_2}{\partial x_1}+g_1\nonumber\\
f_2&&=g_2-u_s\partial_{\bar x_2}v\f{\partial\bar
x_2}{\partial x_1}+\va[2\partial_{\bar x_1\bar x_2}v\f{\partial\bar
x_2}{\partial x_1}+\partial_{\bar x_2}v\f{\partial^2\bar{x}_2}{\partial x_1^2}+\partial_{\bar
x_2}^2v(\f{\partial\bar x_2}{\partial x_2})^2]\nonumber\\
&& +\va(1+\la)\partial_{\bar
x_2}v\f{\partial^2\bar x_2}{\partial x_2^2}+\va(1+\la)[(\f{\partial\bar x_2}{\partial
x_2})^2-1]\partial_{\bar
x_2}^2v+\va\la\partial_{\bar x_1\bar x_2}u(\f{\partial\bar
x_2}{\partial x_2}-1)\nonumber\\
&&+\va\la(u_{\bar x_2\bar x_2}\f{\partial\bar x_2}{\partial x_2}\f{\partial\bar x_2}{\partial x_1}+u_{\bar x_2}\f{\partial^2\bar x_2}{\partial x_1\partial x_2})
-c^2\rho_{\bar x_2}(\f{\partial\bar x_2}{\partial x_2}-1).
\end{eqnarray}


\subsection{Existence of Approximate Solutions  in Hilbert Space}

The linear system (\ref{linear.0}) is a mixed system of hyperbolic-elliptic type, as the velocity $\mau$ satisfies an elliptic system, while the density $\rho$ satisfies a transport equation. We will employ  Larey-Schauder fixed point theory, which can be found in chapter 11 in \cite{G-T},  to prove the existence of solutions to system (\ref{linear.0}). To construct a compact mapping of $\mau$ from $H^1\rightarrow H^1$,  we will first consider the following approximate system:
\begin{eqnarray}\label{app}\begin{cases}
\text{div}\mathbf{u}^\delta+u_s\rho_{x_1}=f_0^\delta\qquad\text{in}~\Omega,\\
u_su_{x_1} +u_{sx_2}v -\va\Delta u
-\va\lambda\partial_{x_1}\mathrm{div}\mathbf{u} +c^2\rho_{x_1} =f_1^\delta\qquad\text{in}~\Omega,\\
u_sv_{x_1} -\va\Delta v
-\va\lambda\partial_{x_2}\mathrm{div}\mathbf{u} +c^2\rho_{x_2} =f_2^\delta\qquad\text{in}~\Omega,
\end{cases}\end{eqnarray}
with boundary condition
\begin{eqnarray}
\rho|_{x_1=0}=u|_{\partial\Omega}=v|_{\partial\Omega}=0,\label{boundary.app}
\end{eqnarray}
here $\mau^\delta, \mathbf{f}^\delta,f_0^\delta$ are the standard mollification of $\mau,\mathbf{f}$ and $f_0$.

First for $t\in [0,1]$,  we consider the momentum equations to be an elliptic system in $\mau$ and the mass equation a transport equation in $\rho$ to have:
\begin{eqnarray}\label{aprt}\begin{cases}
\text{div}\mathbf{u}^\delta+u_s\rho_{x_1}=f_0^\delta\qquad\text{in}~\Omega,\\
 -\va\Delta u -\va\la\partial_{x_1}\dv\mau =t[f_1^\delta-u_su_{x_1}-u_{sx_2}  v-c^2\rho_{x_1}]\ \text{in}\ \Omega,\\
 -\va\Delta v -\va\la\partial_{x_2}\dv\mau =t[f_2^\delta-u_sv_{x_1}-c^2\rho_{x_2}]\ \text{in}\ \Omega,\\
\rho|_{x_1=0}=0,\ \mau |_{\partial\Omega}=0.
\end{cases}\end{eqnarray}

\subsubsection{Apriori Estimates}

\begin{lemma}\label{lemu} Assume that $f_0\in L^2(\Omega),\ \mathbf{f}\in H^{-1}(\Omega)$,  $(\mau,\rho)$ is the strong solution to the approximate system (\ref{aprt}), then
we have the following estimate:
\begin{eqnarray}
&&t\|\mau\|^2+t\|\rho\|^2+\va\|\sqrt{L-x_1}\nabla\mau\|^2\nonumber\\
\leq&& C(\|f_0\|_{L^2}+\va^2\la^2\|v_{x_2}\|^2+\delta\|\mau\|_{H^1_0}^2+ |(t(L-x_1)\mau,\mathbf{f}^\delta)|,\label{l2}
\end{eqnarray}
here the constant $C$ is independent of $\va,\delta,t$.
\end{lemma}
\begin{proof}
We multiply $(\ref{aprt})_2$ with $(L-x_1)u$ and $(\ref{aprt})_3$ with $(L-x_1)v$ to have
\begin{eqnarray*}
&&\int_\Omega[tu_su_{x_1}+tu_{sx_2}v-\va\Delta u-\va\la\partial_{x_1}\mathrm{div}\mathbf{u}+c^2t\rho_{x_1}](L-x_1)udx_1dx_2\nonumber\\
&&+\int_\Omega[tu_sv_{x_1}-\va\Delta v
-\va\la\partial_{x_2}\mathrm{div}\mathbf{u}+tc^2\rho_{x_2}](L-x_1)vdx_1dx_2\nonumber\\
=&&\f12t\int_\Omega(u_su^2+u_sv^2)dx_1dx_2-\f12t\int_\Omega(L-x_1)u_{sx_1}(u^2+v^2)dx_1dx_2\nonumber\\
&&+\va\int_\Omega(L-x_1)[(1+\la)u_{x_1}^2+u_{x_2}^2+v_{x_1}^2+(1+\la)v_{x_2}^2+2\la u_{x_1}v_{x_2}]dx_1dx_2\nonumber\\
&&+c^2t\int_\Omega\rho udx_1dx_2-c^2t\int_\Omega(L-x_1)\rho\dv\mau dx_1dx_2+t\int_\Omega u_{sx_2}vu(L-x_1)dx_1dx_2\nonumber\\
&&-\va\la\int_\Omega uv_{x_2}dx_1dx_2\nonumber\\
=&& (t(L-x_1)\mau,\mathbf{f}^\delta).
\end{eqnarray*}
Using $(\ref{aprt})_1$ we have
\begin{eqnarray*}
&&-c^2t\int_\Omega(L-x_1)\rho\dv\mau dx_1dx_2\nonumber\\
=&&-c^2t\int_\Omega(L-x_1)\rho(f_0^\delta-u_s\rho_{x_1}+\dv\mau-\dv\mau^\delta)dx_1dx_2\nonumber\\
=&&\f12c^2t\int_\Omega [u_s\rho^2-(L-x_1)u_{sx_1}\rho^2]dx_1dx_2-c^2t\int_\Omega(L-x_1)\rho(\dv\mau-\dv\mau^\delta) dx_1dx_2\\
&&+c^2t\int_\Omega(L-x_1)\rho f_0^\delta dx_1dx_2.
\end{eqnarray*}
The supersonic condition implies that
\begin{eqnarray*}
&&\f12t\int_\Omega u_su^2dx_1dx_2+c^2t\int_\Omega\rho udx_1dx_2+\f12c^2t\int_\Omega u_s\rho^2dx_1dx_2\nonumber\\
=&&\f12t\int_\Omega (u_s-c)u^2dx_1dx_2+\f12tc\int_\Omega[u^2+2c\rho u+c^2\rho^2]dx_1dx_2\nonumber\\
&&+\f12c^2t\int_\Omega (u_s-c)\rho^2dx_1dx_2\nonumber\\
=&&\f12t\int_\Omega (u_s-c)u^2dx_1dx_2+\f12tc\int_\Omega(c\rho+u)^2dx_1dx_2+\f12c^2t\int_\Omega (u_s-c)\rho^2dx_1dx_2\nonumber\\
\geq&&\f12t\int_\Omega (u_s-c)u^2dx_1dx_2+\f12c^2t\int_\Omega (u_s-c)\rho^2dx_1dx_2.
\end{eqnarray*}
Combing above, we obtain
\begin{eqnarray*}
&&\f12t\int_\Omega [(u_s-c)u^2+u_sv^2]dx_1dx_2+\f12c^2t\int_\Omega (u_s-c)\rho^2\\
&&+\va\int_\Omega(L-x_1)[\la(\dv\mau)^2+u_{x_1}^2+u_{x_2}^2+v_{x_1}^2+v_{x_2}^2]dx_1dx_2\\
\lesssim&&|-t\int_\Omega u_{sx_2}vu(L-x_1)dx_1dx_2+ (t(L-x_1)\mau,\mathbf{f}^\delta)+c^2t\int_\Omega(L-x_1)\rho f_0^\delta dx_1dx_2\\
&&-c^2t\int_\Omega(L-x_1)\rho(\dv\mau-\dv\mau^\delta) dx_1dx_2+\va\la\int_\Omega uv_{x_2}dx_1dx_2\\
&&+\f12t\int_\Omega(L-x_1)u_{sx_1}(u^2+v^2)dx_1dx_2+\f12c^2t\int_\Omega (L-x_1)u_{sx_1}\rho^2dx_1dx_2|\\
\lesssim&&tL\|\mau\|^2+\|f_0^\delta\|^2+Lt\|\rho\|^2+\va\la\|u\|\|v_{x_2}\|+\delta\|\mau\|_{H^1_0}^2+ |(t(L-x_1)\mau,\mathbf{f}^\delta)|,
\end{eqnarray*}
and estimate (\ref{l2}) follows immediately.
\end{proof}
\begin{lemma}\label{lemen}
 Assume that $f_0\in L^2(\Omega),\ \mathbf{f}\in H^{-1}(\Omega)$,  $(\mau,\rho)$ is the strong solution to the approximate system (\ref{aprt}), then we have the following estimate:
\begin{eqnarray}
t|\rho(L,\cdot)|^2+\va\|\nabla\mau\|^2\leq C(\|f_0^\delta\|^2+t\|\mau\|^2+ t\|\rho\|^2+|(\mau,\mathbf{f}^\delta)|),\label{energy}
\end{eqnarray}
here the constant $C$ is independent of $\va,\delta,t$.
\end{lemma}
\begin{proof}
We multiply  $(\ref{aprt})_2$ with $u$ and  $(\ref{aprt})_3$ with $v$ to have
\begin{eqnarray*}
&&\int_\Omega[tu_su_{x_1}+tu_{sx_2}v-\va\Delta u-\la\va\partial_{x_1}\mathrm{div}\mathbf{u}+c^2t\rho_{x_1}]udx_1dx_2\nonumber\\
&&+\int_\Omega[tu_sv_{x_1}-\va\Delta v
-\va\la\partial_{x_2}\mathrm{div}\mathbf{u}+c^2t\rho_{x_2}]vdx_1dx_2\nonumber\\
=&& -t\int_\Omega u_{sx_1}(u^2+v^2)dx_1dx_2+t\int_\Omega u_{sx_2}vudx_1dx_2+\va\int_\Omega[|D\mau|^2+\la(\dv\mau)^2]dx_1dx_2\\
 &&-c^2t\int_\Omega\rho\dv\mau dx_1dx_2 \nonumber\\
=&&t(\mau,\mathbf{f}^\delta).
\end{eqnarray*}
Using equation $(\ref{aprt})_1$ we have
\begin{eqnarray*}
&&-c^2t\int_\Omega\rho\dv\mau dx_1dx_2=c^2t\int_\Omega\rho(u_s\rho_{x_1}+\dv\mau^\delta-\dv\mau-f_0^\delta) dx_1dx_2\\
=&&\f12c^2t\int_0^2u_s\rho^2|_{x_1=L}dx_2-\f12c^2t\int_\Omega u_{sx_1}\rho^2dx_1dx_2\\
&&+c^2t\int_\Omega\rho(\dv\mau^\delta-\dv\mau-f_0^\delta)dx_1dx_2
\end{eqnarray*}
Combing above we have
\begin{eqnarray*}
&&t\int_0^2 u_s\rho^2|_{x_1=L}dx_2+\va\int_\Omega|D\mau|^2dx_1dx_2\\
\lesssim&&|-t\int_\Omega u_{sx_2}vudx_1dx_2+c^2t\int_\Omega\rho(\dv\mau^\delta-\dv\mau-f_0^\delta)dx_1dx_2\\
&&+t(\mau,\mathbf{f}^\delta)+\f12c^2t\int_\Omega u_{sx_1}\rho^2dx_1dx_2+t\int_\Omega u_{sx_1}(u^2+v^2)dx_1dx_2|\\
\lesssim&&t\|\mau\|^2+t\delta\|\rho\|\|\mau\|_{H^1_0}+t\|\rho\|\|f_0^\delta\|+\va t\|\rho\|^2+|(\mau,\mathbf{f}^\delta)|,
\end{eqnarray*}
and estimate  (\ref{energy}) follows immediately.
\end{proof}

\subsubsection{Existence of Solutions to the Approximate System}
\begin{theorem}\label{thapp}
For given $f_0\in L^2(\Omega),\ \mathbf{f}=(f_1,f_2)\in (H^{-1}(\Omega))^2$,  $0<\delta\ll\va$, there exists a unique solution $\mau\in H^1_0(\Omega)\cap W^{2,p}(\Omega) $ and $\rho\in W^{1,p}(\Omega)$ to system (\ref{app})-(\ref{boundary.app}) with the following estimates:
\begin{equation}\|\mau\|^2+\|\rho\|^2+|\rho(L,\cdot)|^2+\va\|\nabla\mau\|^2\leq C(\va)(\|f_0\|^2+\|\mathbf{f}\|_{H^{-1}(\Omega)}),\label{energyes}\end{equation}
here the constant $C$ is independent of $\delta$.
\end{theorem}

\begin{proof}
First for given $f_0\in L^2(\Omega)$, we define the operator $S:\bar\mau=(\bar u,\bar v)\in H_0^1(\Omega)\rightarrow\rho\in W^{1,p}(\Omega)$ by
$$\rho(x_1,x_2)=\int_0^{x_1}\f1{u_s}[f^\delta_0(s,x_2)-\dv\bar{\mau}^\delta(s,x_2)]ds.$$
Obviously we have $$\text{div}\bar{\mathbf{u}}^\delta+u_s\rho_{x_1}=f^\delta_0\ \text{in}~\Omega,\ \rho\in W^{1,p}(\Omega),\ \rho|_{x_1=0}=0.$$

Then, for $t\in [0,1], \bar\mau=(\bar u,\bar v)\in H_0^1(\Omega)$, we define the mapping $T: H_0^1\times [0,1]\rightarrow H_0^1$ with $T(\bar\mau,t)=\mathbf{u}=(u,v) $  by the solution of the following elliptic system:
\begin{eqnarray}\label{appt}\begin{cases}
-\va\Delta u -\va\la\partial_{x_1}\dv\mau =t[f_1^\delta-u_s\bar u_{x_1} -u_{sx_2}\bar v-c^2S(\bar\mau)_{x_1}]\ \text{in}\ \Omega,\\
 -\va\Delta v -\va\la\partial_{x_2}\dv\mau =t[f_2^\delta-u_s\bar v_{x_1}-c^2S(\bar\mau)_{x_2}]\ \text{in}\ \Omega,\\
\mau |_{\partial\Omega}=0.
\end{cases}\end{eqnarray}
Since  $t[f_1^\delta-u_s\bar u_{x_1}-u_{sx_2}\bar v-c^2S(\bar\mau)_{x_1}],\ t[f_2^\delta-u_s\bar u_{x_1}-c^2S(\bar\mau)_{x_2}]\in L^2(\Omega)$, The Lax-Milgram theorem implies that there exists a unique weak solution $\mau \in H_0^1(\Omega)$ to system (\ref{appt}) and obviously we have $T(\bar\mau,0)=(0,0)$ for any $\bar\mau\in H_0^1(\Omega)$. Besides, by the theory of Lamm$\acute{e}$ system in non-smooth domains(e.g.  Theorem 3.8.1 in \cite{K-M-R}), we have $\mau\in H^2(\Omega)$ and $T$ is a compact mapping.

Finally, considering the fixed point satisfying $T(\mau,t)=\mau$, we have by Lemma \ref{lemu} and Lemma  \ref{lemen}:
\begin{eqnarray}
&&t|\rho(L,\cdot)|^2+t\|\mau\|^2+ t\|\rho\|^2+\va\|\nabla\mau\|^2\nonumber\\
\leq&& C(\|f_0^\delta\|^2+| (t(L-x_1)\mau,\mathbf{f}^\delta)|+|(\mau,\mathbf{f}^\delta)|).\label{L2}
\end{eqnarray}
Consequently we have
\begin{equation}t|\rho(L,\cdot)|^2+t\|\mau\|^2+ t\|\rho\|^2+\va\|\nabla\mau\|^2\leq C(\va)(\|f_0\|^2+\|\mathbf{f}\|_{H^{-1}}),\label{h5}\end{equation}
here the constant $C$ is independent of $t,\delta$. The Leray-Shauder fixed point theorem implies that there exists a unique solution $\mau\in H^1_0(\Omega)\cap W^{2,p}(\Omega)$ to system (\ref{app})-(\ref{boundary.app}). Finally,
(\ref{energyes}) follows by taking $t=1$ in (\ref{h5}).
\end{proof}

\subsection{Estimates of the Approximate Solutions in $W^{2,p}(\Omega)$ }

In this subsection we will study the  $W^{2,p}(\Omega)$ estimates of the approximate solutions.  First by the observation  that the density on the boundary of the domain can be well controlled, we can homogenize the boundary value of the density. Then we  construct  a function $W=(W^1,W^2)$ defined in (\ref{lamme.3}) that satisfies an inhomogeneous elliptic system with homogeneous Dirichlet boundary condition on the boundary of the first quadrant. In fact, $W$ can be split it into two parts: $W_1=(W_{11}, W_{12})$ defined in (\ref{w1})  and $W_2=(W_{21},W_{22})$ defined in (\ref{w2}),(\ref{w}). After homogenizing the boundary value of $\rho$ to a new function $\hat \rho$, we can reduce the expression of $W_1$ to the convolution of $\nabla\hat\rho$ and the fundamental solution of Laplace operator. Then we can construct $W_2$ by use of the  the Green's function of Laplace operator in the first quadrant. Finally the $L^p$ estimates of the approximate solutions follows from the Calderon-Zygmund theory as well as a careful bootstrap argument.
The main result of this section reads as follows.
\begin{theorem}\label{2p}Assume that $f_0\in W^{1,p}(\Omega),\ \mathbf{f}=(f_1,f_2)\in (L^p(\Omega))^2$ in (\ref{app}).  Let $(\mathbf{u},\rho)\in W^{2,p}(\Omega)\times W^{1,p}(\Omega)$ be the solution to the system (\ref{app}) established in Theorem \ref{thapp}, then we have
 the following estimate:
\begin{eqnarray}
\|\rho\|_{1,p;\Omega}+\va\|\mathbf{u}\|_{2,p;\Omega}\leq
C[\|\mathbf{f}\|_{p;\Omega}+\va\|f_0\|_{W^{1,p}}+\va^{-\f32+\f2p-\f12\sigma}(\|\mathbf{f}\|+\|f_0\|)],
\end{eqnarray}
where $\sigma>0$ is any constant small enough, the constant C is independent of
$\va,\delta$.
\end{theorem}
The proof of Theorem \ref{2p} will be split into several Lemmas.

\begin{lemma}\label{rhoy}
Let $(\mathbf{u},\rho)\in W^{2,p}(\Omega)\times W^{1,p}(\Omega)$ be the solution to the system (\ref{app}). Assume that all the assumptions in Theorem \ref{2p} are satisfied, then we have the following estimate:
\begin{eqnarray*}
&&\|\rho_{x_2}\|_{L^p}+\va^{\f1p}|\rho_{x_2}(L,\cdot)|_{L^p}\\
\lesssim&&\va(\|f_0\|_{W^{1,p}}+\|\rho_{x_1}\|_{L^p})+\|f_2\|_{L^p}+\va\|(u_{x_2}-v_{x_1})_{x_1}\|_{L^p}+\|v_{x_1}\|_{L^p}+\va\delta\|\nabla^2\mau\|_{L^p},
\end{eqnarray*}
here the constant $C$ is independent of $\va,\delta$.
\end{lemma}
\begin{proof}
First we differentiate $(\ref{app})_1$ with respect to $x_2$ to have
\begin{equation}\partial_{x_2}\dv\mau^\delta+u_{sx_2}\rho_{x_1}+u_s\rho_{x_1x_2}=f_{0x_2}^\delta,\label{rhoy.1}\end{equation}
which combined with $(\ref{app})_3$ implies
\begin{eqnarray}
&&c^2\rho_{x_2}+\va(1+\la)u_s\rho_{x_1x_2}\nonumber\\
=&&f_2^\delta+\va(1+\la)[\partial_{x_2}\dv\mau-\partial_{x_2}\dv\mau^\delta+f_{0x_2}^\delta-u_{sx_2}\rho_{x_1}]-u_sv_{x_1}-\va(u_{x_2}-v_{x_1})_{x_1}.\nonumber\\\label{deny}
\end{eqnarray}
Then we multiply (\ref{deny}) with $\rho_{x_2}|\rho_{x_2}|^{p-2}$ to have
\begin{eqnarray*}
&&c^2\int_\Omega|\rho_{x_2}|^pdx_1dx_2+\va(1+\la)\int_\Omega u_s\rho_{x_1x_2}\rho_{x_2}|\rho_{x_2}|^{p-2}dx_1dx_2\\
=&&c^2\int_\Omega|\rho_{x_2}|^pdx_1dx_2+\f1p\va(1+\la)\int_0^2u_s|\rho_{x_2}|^p|_{x_1=L}dx_2\\
&&-\f1p\va(1+\la)\int_\Omega u_{sx_1}|\rho_{x_2}|^pdx_1dx_2\\
=&&\int_\Omega[f_2^\delta+\va(1+\la)(\partial_{x_2}\dv\mau-\partial_{x_2}\dv\mau^\delta+f_{0x_2}^\delta-u_{sx_2}\rho_{x_1})]\rho_{x_2}|\rho_{x_2}|^{p-2}dx_1dx_2\\
&&\int_\Omega[-u_sv_{x_1}-\va(u_{x_2}-v_{x_1})_{x_1}]\rho_{x_2}|\rho_{x_2}|^{p-2}dx_1dx_2\\
\lesssim&&\|\rho_{x_2}\|_{L^p}^{p-1}[\va(1+\la)(\|f_0\|_{W^{1,p}}+\|\rho_{x_1}\|_{L^p})+\|f_2\|_{L^p}+\|v_{x_1}\|_{L^p}+\va\|(u_{x_2}-v_{x_1})_{x_1}\|_{L^p}]\\
&&+\va\delta\|\rho_{x_2}\|_{L^p}^{p-1}\|\nabla^2\mau\|_{L^p}.
\end{eqnarray*}
Consequently we have
\begin{eqnarray*}
&&\|\rho_{x_2}\|_{L^p}+\va^{\f1p}|\rho_{x_2}(L,\cdot)|_{L^p}\\
\lesssim&&\va(\|f_0\|_{W^{1,p}}+\|\rho_{x_1}\|_{L^p})+\|f_2\|_{L^p}+\va\|(u_{x_2}-v_{x_1})_{x_1}\|_{L^p}+\|v_{x_1}\|_{L^p}+\va\delta\|\nabla^2\mau\|_{L^p}
\end{eqnarray*}
and we finish the proof of the Lemma.
\end{proof}

From Lemma \ref{rhoy} we find that to prove Theorem \ref{2p}, the key point is to obtain the bound of $\|(u_{x_2}-v_{x_1})_{x_1}\|_{L^p}$. To achieve this, we  write  $(\ref{app})_2$ and $(\ref{app})_3$ into the form of the  Lam$\acute{e}$ system:
\begin{eqnarray}\label{lamme.0}\begin{cases}
 -\va\Delta u
-\va\lambda\partial_{x_1}\mathrm{div}\mathbf{u}  =-c^2\rho_{x_1}+f_1^\delta-u_su_{x_1} -u_{sx_2}v\qquad\text{in}~\Omega,\\
-\va\Delta v
-\va\lambda\partial_{x_2}\mathrm{div}\mathbf{u}  =-c^2\rho_{x_2}+f_2^\delta-u_sv_{x_1} \qquad\text{in}~\Omega.
\end{cases}\end{eqnarray}

We define $\varrho$ by the following elliptic problem:
\begin{equation}\label{homo}\begin{cases}
\Delta \varrho=0,\ \text{in}\ \Omega,\\\varrho=\rho,\ \text{on}\ \partial\Omega.
\end{cases}\end{equation}
The classical theory of elliptic equations implies that there exists a unique solution $\varrho\in W^{1,p}(\Omega)$ to system (\ref{homo}) with the following estimate:
 \begin{equation}\|\varrho\|_{1,p;\Omega}\leq C|\rho|_{1-\f1p,p;\partial\Omega}.\label{rho.1}\end{equation}
 Then using equation $(\ref{app})_1$, Lemma \ref{rhoy} and the trace theory, we have the following estimate for $\varrho$:
\begin{lemma}\label{lp}
Let $(\rho,\mau)\in W^{1,p}(\Omega)\times (W^{2,p}(\Omega))^2$ be the solution to the approximate system (\ref{app}). $\sigma>0$ is a constant sufficiently small, $\varrho\in W^{1,p}(\Omega)$ is defined in (\ref{homo}), then we have the following estimate:
\begin{eqnarray}
\|\varrho\|_{W^{1,p}}\leq&&\va^{1+\f\sigma{10}}[\|f_0\|_{w^{1,p}}+\|\rho_x\|_{L^p}+\|\mau\|_{W^{2,p}}]+\|f_2\|_{L^p}\nonumber\\
&&+\va^{-\f32+\f2p-\f12\sigma}(\|\mathbf{f}\|+\|f_0\|)+\va\delta\|\nabla^2\mau\|_{L^p}\label{p.3}
\end{eqnarray}
\end{lemma}
\begin{proof}
As $\rho|_{x_1=0}=0$, we have by (\ref{rho.1}):
\begin{eqnarray}
\|\varrho\|_{W^{1,p}}\lesssim |\rho(L,\cdot)|_{W^{1-1/p,p}(0,2)}+|\rho(\cdot,0)|_{W^{1-1/p,p}(0,L)}+|\rho(\cdot,2)|_{W^{1-1/p,p}(0,L)}\label{rho.2}.
\end{eqnarray}
In the following we will estimate the right side of (\ref{rho.2}) term by term.
First using (\ref{L2}), Lemma \ref{rhoy},  the extension theory, the Young's inequality and Gagliardo-Nirenberg inequality  we have
\begin{eqnarray*}
&&|\rho(L,\cdot)|_{W^{1-1/p,p}(0,2)}\lesssim|\rho_{x_2}(L,\cdot)|^{1-\f2{3p-2}}_{L^p}|\rho(L,\cdot)|^{\f2{3p-2}}_{L^2}+|\rho(L,\cdot)|_{L^2}\\
\lesssim&& \sigma_1|\rho_{x_2}(L,\cdot)|_{L^p}+\sigma_1^{2-\f32p}|\rho(L,\cdot)|_{L^2}\nonumber \\
\lesssim&&\sigma_1\va^{-\f1p}[\va(\|f_0\|_{W^{1,p}}+\|\rho_{x_1}\|_{L^p})+\|f_2\|_{L^p}+\va\|(u_{x_2}-v_{x_1})_{x_1}\|_{L^p}+\|v_{x_1}\|_{L^p}\nonumber\\
&&+\va\delta\|\nabla^2\mau\|_{L^p}]+\sigma_1^{2-\f32p}(\|\mathbf{f}\|+\|f_0\|).
\end{eqnarray*}
For  $\sigma>0$ small enough and $M>0$ big enough, we have by taking $\sigma_1=\va^{\f1p+\f\sigma M}$:
\begin{eqnarray}
&&|\rho(L,\cdot)|_{W^{1-1/p,p}(0,2)}\nonumber \\
\lesssim&&\va^{\f\sigma M}[\va(\|f_0\|_{w^{1,p}}+\|\rho_{x_1}\|_{L^p})+\|f_2\|_{L^p}+\va\|(u_{x_2}-v_{x_1})_{x_1}\|_{L^p}+\|v_{x_1}\|_{L^p}\nonumber\\
&&+\va\delta\|\nabla^2\mau\|_{L^p}]+\va^{\f2p-\f32+\f\sigma{2M}-\f3{8M}p\sigma}(\|\mathbf{f}\|+\|f_0\|)\nonumber\\
\lesssim&&\va^{1+\f\sigma M}[\|f_0\|_{W^{1,p}}+\|\rho_{x_1}\|_{L^p}+\|\nabla\mau_{x_1}\|_{L^p}]+\|f_2\|_{L^p}+\va^{\f2p-\f32-\f1{2M}\sigma}(\|\mathbf{f}\|+\|f_0\|)\nonumber\\
&&+\va\delta\|\nabla^2\mau\|_{L^p},
\end{eqnarray}
here we have used the inequality  from  Gagliardo-Nirenberg inequality  that for $2<p<\infty,\ \omega\in W^{1,p}(\Omega)$ we have
 \begin{eqnarray}
 \|\omega\|_{p;\Omega}\lesssim\|\nabla\omega\|_{L^p(\Omega)}^{\f{p-2}{2(p-1)}}\|\omega\|_{L^2(\Omega)}^{\f{p}{2(p-1)}}+\|\omega\|_{2;\Omega}
 \leq \va^{1+\f\sigma M}\|\nabla\omega\|_{p;\Omega}+\va^{\f{2}{p}-1-\f\sigma M}\|\omega\|_{L^2}.\label{du}
 \end{eqnarray}
Similarly as above we have
$$|\rho(\cdot,2)|_{W^{1-1/p,p}(0,L)}\leq \|\rho_{x_1}(\cdot,2)\|_{L^p}^{a}\|\rho(\cdot,2)\|_{L^\infty}^{1-a}+\|\rho(\cdot,2)\|_{L^p},$$
where $a=1-\f1{p-1}$.  By equation $(\ref{app})_1$, the trace theory and Gagliardo-Nirenberg inequality we have
\begin{eqnarray*}&&\|\rho(\cdot,2)\|_{L^\infty}=\sup_{x_1\in [0,L]}|\int_0^{x_1}\rho_{x_1}(s,2)ds|
\lesssim|\rho_{x_1}(\cdot,2)|_{L^1(0,L)}\lesssim |(f_0-\dv\mau)(\cdot,2)|_{L^1}\\
\lesssim&& \|f_0\|_{W^{1-\sigma,1+2\sigma}(\Omega)}+\|\dv\mau\|_{W^{1-\sigma,1+2\sigma}(\Omega)}\\
\lesssim&&\|\dv\mau\|_{W^{1,p}(\Omega)}^{a_1}\|\dv\mau\|_{L^2(\Omega)}^{1-a_1}+\|f_0\|_{W^{1,p}}^{a_1}\|f_0\|_{L^2}^{1-a_1}+\|\dv\mau\|_{L^2(\Omega)}+\|f_0\|_{L^2(\Omega)}
\end{eqnarray*}
where $a_1=\f{(3-2\sigma)p\sigma}{2(1+2\sigma)(p-1)}$. In the same way we can obtain
\begin{eqnarray*}&&\|\rho_{x_1}(\cdot,2)\|_{L^p}= |\f1{u_s}(f_0^\delta-\dv\mau^\delta)(\cdot,2)|_{L^p}\lesssim \|f_0\|_{W^{\f1p+\sigma,p}(\Omega)}+\|\dv\mau\|_{W^{\f1p+\sigma,p}(\Omega)}\\
\lesssim&& \|\dv\mau\|_{W^{1,p}(\Omega)}^{a_2}\|\dv\mau\|_{L^2(\Omega)}^{1-a_2}+\|f_0\|_{W^{1,p}}^{a_2}\|f_0\|_{L^2}^{1-a_2}+\|\dv\mau\|_{L^2(\Omega)}+\|f_0\|_{L^2(\Omega)},
\end{eqnarray*}
where $a_2=\f12+\f{p\sigma}{2(p-1)}$.
Combing above and using Young's inequality, we obtain that
\begin{eqnarray}
&&|\rho(\cdot,2)|_{W^{1-1/p,p}(0,L)}\nonumber\\
\leq &&\|\dv\mau\|_{W^{1,p}(\Omega)}^{a_2a+a_1(1-a)}\|\dv\mau\|_{L^2(\Omega)}^{a(1-a_2)+(1-a)(1-a_1)}+\|\dv\mau\|_{L^2(\Omega)}+\|f_0\|_{L^2(\Omega)}\nonumber\\
&&+\|f_0\|_{W^{1,p}(\Omega)}^{a_2a+a_1(1-a)}\|f_0\|_{L^2(\Omega)}^{a(1-a_2)+(1-a)(1-a_1)}\nonumber\\
\leq&& \sigma_2\|\dv\mau\|_{W^{1,p}(\Omega)}+\sigma_2^{-a_3}\|\dv\mau\|_{L^2(\Omega)}+\sigma_2\|f_0\|_{W^{1,p}(\Omega)}+\sigma_2^{-a_3}\|f_0\|_{L^2(\Omega)}
\end{eqnarray}
where direct computation shows that
 $$a_3=\f{a_2a+a_1(1-a)}{a(1-a_2)+(1-a)(1-a_1)}=\f{(a_2-a_1)a+a_1}{a(a_1-a_2)+1-a_1}\leq \f{p-2}{p}+5\sigma.$$
If we take $\sigma_2=\va^{1+\f\sigma{10}}$, then we have
\begin{eqnarray}
&&|\rho(\cdot,2)|_{W^{1-1/p,p}(0,L)}\nonumber\\
\leq&& \va^{1+\f\sigma{10}}\|\dv\mau\|_{W^{1,p}(\Omega)}+\va^{-(\f{p-2}{p}+\f\sigma2)}\|\dv\mau\|_{L^2(\Omega)}+\va^{1+\f\sigma{10}}\|f_0\|_{W^{1,p}(\Omega)}\nonumber\\
&&+\va^{-(\f{p-2}{p}+\f\sigma2)}\|f_0\|_{L^2(\Omega)}\nonumber\\
\leq&& \va^{1+\f\sigma{10}}(\|\dv\mau\|_{W^{1,p}(\Omega)}+\|f_0\|_{W^{1,p}(\Omega)})+\va^{-\f32+\f2p-\f12\sigma}(\|f_0\|_{L^2(\Omega)}+\|\mathbf{f}\|_{L^2})
\end{eqnarray}
The estimates for $|\rho(\cdot,0)|_{W^{1-1/p,p}(0,L)}$ can be obtained in the same way and we have proved estimate (\ref{p.3}).
\end{proof}
If we denote by $\tilde\rho=\rho-\varrho$, then we have $\tilde\rho\in W^{1,p}_0(\Omega)$ and system (\ref{lamme.0}) can be written into the following form:
\begin{eqnarray}\label{lamme.4}\begin{cases}
 -\va\Delta u
-\va\lambda\partial_{x_1}\mathrm{div}\mathbf{u}  =-c^2\tilde\rho_{x_1}+\tilde f_1-u_su_{x_1} -u_{sx_2}v\qquad\text{in}~\Omega,\\
-\va\Delta v
-\va\lambda\partial_{x_2}\mathrm{div}\mathbf{u}  =-c^2\tilde\rho_{x_2}+\tilde f_2-u_sv_{x_1} \qquad\text{in}~\Omega,
\end{cases}\end{eqnarray}
here $\tilde{\mathbf{f}}=(\tilde f_1,\tilde f_2)\in L^p(\Omega)$ with $\tilde f_1=f_1^\delta-c^2\varrho_{x_1},\ \tilde f_2=f_2^\delta-c^2\varrho_{x_2}$.

Next, we will study system (\ref{lamme.4}) in the first quadrant. Here we present some known results on the elliptic system which will be used in the proof of Theorem \ref{2p}. For more details we refer to \cite{ADN2}.
First, for $x=(x_1,x_2)\in \mathbb{R}^2$, we consider the following Lamm$\acute{e}$ system
\begin{eqnarray}
\sum^2_{j=1}l_{ij}(\partial)u_j(x)= \va\Delta
u_i+ \lambda\va\partial_i\mathrm{div} u&&=:\phi_i,\quad
i=1,2,\qquad\text{in}~\mathbb{R}^2,\label{lamme.1}
\end{eqnarray}
here $l_{ij}(\partial)$ is a polynomial in $\partial$. Let $L^{jk}$ denote the adjoint to the matrix $\{l_{ij}\}$, i.e.
$$\sum_{j=1}^{N}l_{ij}(\xi)L^{jk}(\xi)=\delta^k_iL(\xi),\qquad i,k=1,...,N,$$
where $L(\xi)=\mathrm{det}\{l_{ij}(\xi )\}= \va^2(1+\la)|\xi|^4$ and $\delta^k_i$ is Kronecker's delta. Then direct computation shows that
\begin{eqnarray}L^{jk}(\xi)=\left(
  \begin{array}{cc}
    L^{11} & L^{12} \\
    L^{21} & L^{22} \\
  \end{array}
\right)=\left(
  \begin{array}{cc}
     \va\xi^2_1+ \va(1+\lambda)\xi_2^2 & - \va\lambda\xi_1\xi_2 \\
   - \va\lambda\xi_1\xi_2 & \va (1+\lambda)\xi_1^2+ \xi_2^2. \\
  \end{array}
\right)\label{L}\end{eqnarray}

For $\phi_i\in C^\infty_0({\mathbb R}^2)$, the functions $u_i$, defined by
\begin{eqnarray}
u_i(x)=\sum^N_{j=1}\int_{{\mathbb R}^2}[L^{ji}(\partial_y)\Gamma(x-y)]\phi_j(y)dy,\qquad\quad
i=1,2,\label{lamme.2}
\end{eqnarray}
satisfy the differential equations (\ref{lamme.1}) in $\mathbb{R}^2$, where
$\Gamma(x)=\f1{8\pi}|x|^2\ln|x|$ is the fundamental solution of the
biharmonic operator $\Delta^2$ in $\mathbb{R}^2$.

 Letting
\begin{eqnarray*}
&&Q_{R}(0)=\{(x_1,x_2)\in \mathbb{R}^2\ |\ |x|<R,\ x_1>0,x_2>0\},\\
&&Q^+=\{(x_1,x_2)\in \mathbb{R}^2\ |\ \ x_1>0,x_2>0\},\\
&&T_1=\{(x_1,0)\in \mathbb{R}^2\ |\ \ x_1>0\}, \\
&&T_2=\{(0,x_2)\in \mathbb{R}^2\ |\ \ x_2>0\}.
\end{eqnarray*}
For given $\hat\rho\in W^{1,p}_0(Q_{2R}(0))$,
 we define $W_1=(W_{11},W_{12})$ with
\begin{eqnarray}
\va W_{1i}(x)=\int_{Q_{2R}(0)}L^{ij}(\partial_y)\Gamma(x-y)\partial_j\hat\rho(y)dy\qquad\text{in}~Q_{2R}(0),~i=1,2,\label{w1}
\end{eqnarray}
here  $L^{ij}$ are defined in (\ref{L}). By (\ref{lamme.1}) and (\ref{lamme.2}), we see that
\begin{equation}
\va\Delta W_1+\va\lambda\nabla\mathrm{div} W_1=\nabla\hat\rho~\text{in}~Q_{2R}(0)\label{w11}.
\end{equation}
\begin{lemma}\label{ep}
Assume that  $\hat\rho\in W^{1,p}_0(Q_{2R}(0))$,  $W_1=(W_{11},W_{12})$ is defined in (\ref{w1}), then we have
$$\va\Delta W_{1}=\nabla\hat\rho\  \text{in}\  Q_{2R}(0),$$
 and
\begin{eqnarray}
&&\va\|W_{11}\|_{2,p;Q_{2R}(0)}+\va\|\partial_1W_{12}\|_{1,p;Q_{2R}(0)}\leq C\|\partial_1\hat\rho\|_{p;Q_{2R}(0)},
\end{eqnarray}
where the constant  $C$ is independent of $\va$.
\end{lemma}
\begin{proof}
 First recalling the definition of $L^{ij}$ in (\ref{L}), since $\hat\rho\in W^{1,p}_0(Q_{2R}(0))$, we have by integration by parts:
\begin{eqnarray}
&&\va W_{11}(x)=\int_{Q_{2R}(0)}L^{11}(\partial_y)\Gamma(x-y)\partial_1\hat\rho(y)dy
+\int_{Q_{2R}(0)}L^{12}(\partial_y)\Gamma(x-y)\partial_2\hat\rho(y)dy\nonumber\\
=&& \int_{Q_{2R}(0)}[\Delta_y\Gamma(x-y)+\la\partial_{y_2y_2}\Gamma(x-y)]\partial_1\hat\rho(y)dy\nonumber\\
&&- \lambda\int_{Q_{2R}(0)}\partial_{2}\hat\rho(y)\partial_{y_1y_2}\Gamma(x-y)dy\nonumber\\
=&& \int_{Q_{2R}(0)}[\Delta_y\Gamma(x-y)+\la\partial_{y_2y_2}\Gamma(x-y)]\partial_1\hat\rho(y)dy\nonumber\\
&&- \lambda\int_{Q_{2R}(0)}\partial_{1}\hat\rho(y)\partial_{y_2y_2}\Gamma(x-y)dy\nonumber\\
=&&\f1{2\pi} \int_{Q_{2R}(0)}\partial_{1}\hat\rho(y)(\ln|x-y|+1)dy\nonumber\\
=&&\f1{2\pi} \int_{Q_{2R}(0)}\partial_{1}\hat\rho(y)\ln|x-y|dy,\label{h1}
\end{eqnarray}
and
\begin{eqnarray}
&&\va W_{12}(x)=\int_{Q_{2R}(0)}L^{21}(\partial_y)\Gamma(x-y)\partial_1\hat\rho(y)dy
+\int_{Q_{2R}(0)}L^{22}(\partial_y)\Gamma(x-y)\partial_2\hat\rho(y)dy\nonumber\\
=&&- \lambda\int_{Q_{2R}(0)}\partial_{y_1y_2}\Gamma(x-y)\partial_1\hat\rho(y)dy+ \int_{Q_{2R}(0)}\Delta\Gamma(x-y)\partial_2\hat\rho(y)dy\nonumber\\
&&+ \la\int_{Q_{2R}(0)}\partial_{y_1y_1}\Gamma(x-y)\partial_2\hat\rho(y)dy\nonumber\\
=&&\f1{2\pi} \int_{Q_{2R}(0)}\partial_{2}\hat\rho(y)(\ln|x-y|+1)dy\nonumber\\
=&&\f1{2\pi} \int_{Q_{2R}(0)}\partial_{2}\hat\rho(y)\ln|x-y|dy.\label{h2}
\end{eqnarray}
Then we have
 $$\va\Delta W_{11}=\partial_1\hat\rho,\   \va\Delta W_{12}=\partial_2\hat\rho \ \ \text{in} \  Q_{2R}(0),$$
  and the Calderon-Zygmund theory implies
 \begin{equation}
\va \|W_{11}\|_{2,p;\Omega}\leq C\|\partial_1\hat\rho\|_{L^p(\Omega)}.
 \end{equation}

Next direct computation shows that:
\begin{eqnarray*}
&&\va\partial_1W_{12}(x)=\f1{2\pi} \int_{Q_{2R}(0)}\partial_2\hat\rho(y)\f{x_1-y_1}{|x-y|^2}dy  \nonumber\\
=&&-\f1{2\pi}\lim_{\sigma\rightarrow 0} \int_{\partial B_{\sigma}(x)}\hat\rho(y)\f{x_1-y_1}{|x-y|^2}\f{x_2-y_2}{|x-y|}ds\nonumber\\
&&-\f1{2\pi}\lim_{\sigma\rightarrow 0} \int_{Q_{2R}(0)/B_{\sigma}(x)}\hat\rho(y)\partial_{12}\ln|x-y|dy\nonumber\\
=&&-\f1{2\pi}\lim_{\sigma\rightarrow 0} \int_0^{2\pi}\f12\hat\rho(x_1+\sigma\cos\theta,x_2+\sigma\sin\theta d\theta)\sin2\theta d\theta\\
&&-\f1{2\pi}\lim_{\sigma\rightarrow 0} \int_{Q_{2R}(0)/B_{\sigma}(x)}\hat\rho(y)\partial_{12}\ln|x-y|dy\nonumber\\
=&&-\f1{2\pi}\lim_{\sigma\rightarrow 0} \int_{Q_{2R}(0)/B_{\sigma}(x)}\hat\rho(y)\partial_{12}\ln|x-y|dy\nonumber\\
=&&\f1{2\pi}\lim_{\sigma\rightarrow 0} \int_{\partial B_{\sigma}(x)}\hat\rho(y)\f{x_1-y_1}{|x-y|}\f{x_2-y_2}{|x-y|^2}ds+ \f1{2\pi}\int_{Q_{2R}(0)}\partial_1\hat\rho(y)\partial_{2}\ln|x-y|dy\nonumber\\
=&& \f1{2\pi}\int_{Q_{2R}(0)}\partial_1\hat\rho(y)\partial_{2}\ln|x-y|dy.
\end{eqnarray*}
Then the Calderon-Zygmund theory implies that
$$\va\|\partial_1W_{12}\|_{1,p;Q_{2R}(0)}\leq
C\|\partial_1\hat\rho\|_{p;Q_{2R}(0)},$$
and we have finished the proof of the lemma.
\end{proof}
%

Based on $W_1$ defined in (\ref{w1}), we will construct a function $W=(W^1,W^2)$  satisfying the following system with homogenous boundary value on $T_1,T_2$:
 \begin{eqnarray}\label{lamme.3}\begin{cases}\va\Delta
W+\lambda\va\nabla\mathrm{div}
W=\nabla\hat\rho+\mathbf{g}\quad\text{in}~Q_{2R}(0),\\
W=0\qquad \text{on}\; T_1\cup T_2.\end{cases}\end{eqnarray}

Step I: Construction of $W^1$:

As we have $\hat\rho\in W^{1,p}_0(Q_{2R}(0))$, if we define $W^1$ by
  the following elliptic  boundary problem:
\begin{equation}\begin{cases}
\va\Delta W^1=\partial_{x_1}\hat\rho,\ \text{in}\ Q_{2R}(0),\\
W^1=0,\ \text{on}\ \partial Q_{2R}(0),
\end{cases}
\end{equation}
then by the theory of elliptic problems we have
$$\va\|W^1\|_{W^{2,p}(Q_{2R}(0))}\leq C \|\hat\rho_{x_1}\|_{L^p(Q_{2R}(0))}.$$
Denoting
\begin{equation}W_{21}=W^1-W_{11},\label{w2}\end{equation}
then $\Delta W_{21}=0$ in $Q_{2R}(0)$, and by Lemma \ref{ep} we have
\begin{equation}\|W_{21}\|_{W^{2,p}(Q_{2R}(0))}\leq\|W_{11}\|_{W^{2,p}(Q_{2R}(0))}+\|W^1\|_{W^{2,p}(Q_{2R}(0))}\leq C \|\partial_{x_1}\hat\rho\|_{L^p(Q_{2R}(0))}.\label{w21}\end{equation}

Step II: Construction of $W^2$:

For $x=(x_1,x_2)\in Q_{2R}(0)$, we denote by
$$x^*=(x_1,-x_2),\ x^{**}=(-x_1,-x_2),\ x^{***}=(-x_1,x_2),$$
and denote by $G(x,y)$ the (Dirichlet) Green's function  for Laplace operator in the first quadrant $Q^+$, i.e.
\begin{eqnarray*}
G(x,y)=\ln|x-y|-\ln|x^*-y|+\ln|x^{**}-y|-\ln|x^{***}-y|,\ x,y\in Q^+.
\end{eqnarray*}
In fact, it is easy to check that  $$\Delta_xG(x,y)=\delta(x-y)\ \text{ for}\ x,y\in Q^+,$$ and for any $x\in Q^+$, when $y=(y_1,0)$ with $y_1>0$, it holds that
\begin{equation}|x-y|=|x^*-y|,\ |x^{**}-y|=|x^{***}-y|,\label{t1}\end{equation}
while when $y=(0,y_2)$ with $y_2>0$, it holds that
\begin{equation}|x-y|=|x^{***}-y|,\ |x^{*}-y|=|x^{**}-y|,\label{t2}\end{equation}
Consequently we have
\begin{equation}G(x,y)=0, \ \forall y\in T_1\cup T_2,\ x\in Q^+.\label{green}\end{equation}
Recalling that $$\va W_{12}(x)=\f1{2\pi} \int_{Q_{2R}(0)}\partial_{2}\hat\rho(y)\ln|x-y|dy,\  \ x\in Q_{2R}(0),$$
if we define $W_{22}$ by
\begin{eqnarray*}
\va W_{22}(x)=\f1{2\pi} \int_{Q_{2R}(0)}\partial_2\hat\rho(y)[- \ln|x^*-y|+ \ln|x^{**}-y|- \ln|x^{***}-y| ]dy.
\end{eqnarray*}
Then for any $x\in Q^+$,
\begin{eqnarray*}
\va\Delta W_{22}(x)=&&\f1{2\pi}\int_{Q_{2R}(0)}\partial_2\hat\rho(y)\Delta_x[- \ln|x^*-y|+ \ln|x^{**}-y|- \ln|x^{***}-y|]dy\\
=&&0.
\end{eqnarray*}
and
\begin{eqnarray*}
&&\va W_{12}(x)+\va W_{22}(x)\\
=&&\f1{2\pi} \int_{Q_{2R}(0)}\partial_2\hat\rho(y)[\ln|x-y|-\ln|x^*-y|+\ln|x^{**}-y|-\ln|x^{***}-y|]dy\\
=&&\f1{2\pi} \int_{Q_{2R}(0)}\partial_2\hat\rho(y)G(x,y)dy\\
=&&0,\  \text{on}\ T_1\cup T_2.\end{eqnarray*}
Finally  we define
\begin{equation}\ W^2=W_{12}+W_{22},\ W_2=(W_{21},W_{22}).\label{w}\end{equation}
Recall the definition of $W_1$ in (\ref{w11}), it is easy to check that $W=(W^1,W^2)=(W_{11}+W_{21},W_{12}+W_{22})=W_1+W_2$ satisfy system (\ref{lamme.3}) with $\mathbf{g}=(g_1,g_2)$ and
\begin{equation}\begin{cases}g_1=\la\va\partial_1\dv W_2=\la\va[\partial_{11}W_{21}+\partial_{21}W_{22}],\\ g_2=\va\la\partial_2\dv W_2=\la\va[\partial_{12}W_{21}+\partial_{22}W_{22}].\label{g}\end{cases}\end{equation}
\begin{lemma}\label{lw}
Assume that  $\hat\rho\in W^{1,p}_0(Q_{2R}(0))$,  $W=(W^1,W^2), \mathbf{g}=(g_1,g_2)$ are defined above, then
we have $W_{22}\in W^{2,p}(\Omega),\ \mathbf{g}\in L^p(\Omega)$ and
\begin{eqnarray}
&&\va\|W_{22}\|_{2,p;Q_{2R}(0)}+\|\mathbf{g}\|_{p;Q_{2R}(0)}\leq C\|\partial_1\hat\rho\|_{p;Q_{2R}(0)},\label{wp}
\end{eqnarray}
where the constant $C$ is independent of $\va$.
\end{lemma}
\begin{proof}

First, by integration by parts we have for any $x=(x_1,x_2)\in Q_{2R}(0)$,
\begin{eqnarray*}
\va\partial_{x_1}W_{22}(x)=&&\f1{2\pi} \int_{Q_{2R}(0)}\partial_{x_1}[-\ln|x^*-y|+\ln|x^{**}-y|-\ln|x^{***}-y|]\partial_2\hat \rho(y)dy\\
=&&\f1{2\pi} \int_{Q_{2R}(0)}\partial_{y_1}[\ln|x^*-y|+\ln|x^{**}-y|-\ln|x^{***}-y|]\partial_2\hat \rho(y)dy\\
=&&\f1{2\pi} \int_{Q_{2R}(0)}\partial_{y_2}[\ln|x^*-y|+\ln|x^{**}-y|-\ln|x^{***}-y|]\partial_{y_1}\hat \rho(y)dy\\
=&&\f1{2\pi} \int_{R^2_+}(\partial_{y_2}\ln|x^*-y|)\partial_1(E\hat \rho)dy-\int_{Q^+}(\partial_{y_2}\ln|x^{***}-y|)\partial_1(E\hat \rho)dy\\
&&+\f1{2\pi} \int_{Q_{2R}(0)}(\partial_{y_2}\ln|x^{**}-y|)\partial_1\hat \rho(y)dy\\
\triangleq&&h_1+h_2+h_3,
\end{eqnarray*}
here we denote by $E\hat \rho\in W_0^{1,p}(\mathbb{R}^2)$ the zero extension of $\hat \rho$ to $\mathbb{R}^2$.
Recall that $x^*=(x_1,-x_2)$, it is easy to check that
\begin{eqnarray*}\Delta h_1=0, \ \text{in}\ \mathbb{R}^2_+, \end{eqnarray*}
and on $ \partial B(0,3R)$ we have
\begin{eqnarray*}|h_1(x)|_{L^\infty}\lesssim \f1R\int_{Q_{2R}(0)}|\partial_1\hat\rho|dy\leq \|\partial_1\hat\rho\|_{L^p(Q_{2R}(0))}.\end{eqnarray*}
Then direct computation shows
\begin{eqnarray}
 h_1(x_1,0)=&&\f1{2\pi} \int_{R^2_+}[\partial_{y_2}\ln|x^*-y|\partial_{y_1}(E\hat\rho)]|_{x_2=0}dy\nonumber\\
 =&&\f1{2\pi} \int_{R^2_+}\f{y_2}{(x_1-y_1)^2+y_2^2}\partial_{y_1}(E\hat\rho)dy\nonumber\\
 =&&\f1{2\pi} [\int_{R^2}\partial_{y_2}\ln|x-y|\partial_{y_1}(E\hat\rho)dy]\big|_{x_2=0},\nonumber\\
  =&&\f1{2\pi} [-\partial_{x_2}\int_{R^2}\ln|x-y|\partial_{y_1}(E\hat\rho)dy]\big|_{x_2=0},
\end{eqnarray}
while for $x\in B(0,3R)$,  the potential function
$$h_0(x)=\f1{2\pi}\int_{R^2}\ln|x-y|\partial_1(E\hat\rho)dy$$
satisfying
$$\|h_0\|_{W^{2,p}(B(0,3R))}\leq \|\partial_1(E\rho)\|_{L^p}\leq \|\partial_1\hat\rho\|_{L^p(Q_{2R}(0))}.$$
Since $h_1$ satisfies the following problem
$$\begin{cases}\Delta h_1=0 \ \text{in}\ \mathbb{R}^2_+,\\h_1|_{x_2=0}=-\f1{2\pi}\partial_{x_2}h_0(x_1,0),\\|h_1(x)|_{L^\infty}\lesssim \|\partial_1\hat\rho\|_{L^p(Q_{2R}(0))},\ \text{on}\  \partial B(0,3R),\end{cases}$$
the  elliptic theory implies that
$$\|h_1\|_{1,p;Q_{2R}}\leq C\|\partial_1\hat\rho\|_{p;Q_{2R}(0)}.$$
Similarly, we check  that $\Delta h_2=0$, $|h_2(x)|_{L^\infty}\lesssim  \|\partial_1\hat\rho\|_{L^p(Q_{2R}(0))},\ \text{on}\  \partial B(0,3R)$
\begin{eqnarray}
 h_2(x)|_{x_1=0}=\f1{2\pi}[\int_{R^2}\partial_{y_2}\ln|x-y|\partial_1(E\rho)dy]|_{x_1=0},
\end{eqnarray}
and we have
$$\|h_2\|_{1,p;Q_{2R}}\leq C\|\partial_1\hat\rho\|_{p;Q_{2R}(0)}.$$
Finally, $h_3$ satisfies
\begin{equation}\begin{cases}
\Delta h_3=0,\ \text{in}\ Q^+,\\
h_3|_{x_1=0}=h_1(0,x_2),\ h_3|_{x_2=0}=h_2(x_1,0), \\
|h_3(x)|_{L^\infty}\lesssim  \|\partial_1\hat\rho\|_{L^p(Q_{2R}(0))},\ \text{on}\  \partial B(0,3R)
\end{cases}\end{equation}
and we have the compatibility condition on the corner $(0,0)$:
$$h_2(0,0)=h_1(0,0)=\f1{2\pi} \int_{Q_{2R}(0)}\f{y_2}{y_1^2+y_2^2}\partial_1\hat \rho(y)dy.$$
The elliptic theory implies that
$$\|h_3\|_{1,p;Q_{2R}(0)}\leq C\|\partial_1\hat\rho\|_{p;Q_{2R}(0)}.$$
Combing above, we have proved that
$$\va\|\partial_1W_{22}\|_{1,p;Q_{2R}(0)}\leq C\|\partial_1\hat\rho\|_{p;Q_{2R}(0)}.$$
Recall that  $\Delta W_{22}=0$ in $Q^+$,  we have
$$\|\partial_{22}W_{22}\|_{L^p;Q_{2R}(0)}\leq \|\partial_{11}W_{22}\|_{L^p;Q_{2R}(0)}\leq C\|\partial_1\hat\rho\|_{p;Q_{2R}(0)}.$$
Finally, the $L^p$ estimate of $\mathbf{g}$ follows immediately from (\ref{w21}),(\ref{g}) and the estimate of $W_{22}$ above.

\end{proof}

$\mathbf{Proof\ of\  Theorem\ \ref{2p}}:$

Letting $\chi(t)$ be the cut-off function defined in  (\ref{cutoff}) and we denote by
\begin{equation}\bar\mau=(\bar u,\bar v)=\chi(\f{x_1}{L})\chi(\f{x_2}2)(u,v),\ \bar\rho=\chi(\f{x_1}{L})\chi(\f{x_2}2)\tilde\rho,\label{h4}\end{equation}
then from system (\ref{lamme.4}) we have

\begin{align*}
\va\Delta\bar{ u}
+\va\lambda\partial_{x_1}\mathrm{div}\bar{\mathbf{u}}&=c^2\bar\rho_{x_1}-\bar f_1,& \text{in}~Q_2(0), \\
\va\Delta\bar v
+\va\lambda\partial_{x_2}\mathrm{div}\bar{\mathbf{u}}&=c^2\bar{\rho}_{x_2}-\bar f_2,& \text{in}~Q_2(0),\\
\bar{\mau}&=0,&\text{on}\ \partial Q_2(0),
\end{align*}
here $\bar\rho\in W_0^{1,p}(Q_2(0))$ and $\bar{\mathbf{f}}=(\bar f_1,\bar f_2)$ with
\begin{align*}
\bar f_1=&\chi(\f{x_1}{L})\chi(\f{x_2}2)(f_1^\delta-u_su_{x_1}-u_{sx_2}v)+\f{c^2}L\chi'(\f{x_1}{L})\chi(\f{x_2}2)\tilde\rho-\va(\f{\chi''}{4}u+\f12\chi'u_{x_2})\chi(\f{x_1}L)\\
&-\va(1+\lambda)(\f{\chi''}{L^2}u+\f2L\chi'u_{x_1})\chi(\f{x_2}2)-\va\lambda\f{\chi'}L(\chi'(\f {x_2}2) v)_{x_2}\\
\bar f_2=&\chi(\f{x_1}{L})\chi(\f{x_2}2)(f_2^\delta-u_sv_{x_1})+\f{c^2}2\chi(\f{x_1}{L})\chi'(\f{x_2}2)\tilde\rho-\va(\f{\chi''}{L^2}v+\f2L\chi'v_{x_1})\chi(\f{x_2}2)\\
&-\va(1+\lambda)(\f{\chi''}{4}v+\f12\chi'v_{x_2})\chi(\f{x_1}L)-\va\lambda\f{\chi'}L(\chi'(\f {x_2}2) u)_{x_2}
\end{align*}
Now let $  W=(  W^1,  W^2)$ with $  W^1=  W_{11}+  W_{12}$ and $  W^2=  W_{21}+  W_{22}$ be as in Lemma \ref{lw} with $\hat\rho$ replaced by $c^2\bar\rho$, then we have
 \begin{eqnarray}
 \begin{cases}
 \va\Delta  W+\lambda\va\nabla\mathrm{div}
W=c^2\nabla \bar\rho+ \mathbf{g} \quad\text{in}~Q_{2 }(0),\\
  W=0\qquad \text{on}\; T_1\cup T_2.\end{cases}\end{eqnarray}
with $ \mathbf{g} =(  g_1,  g_2)$ defined in (\ref{g}).
If we define $\hat \mau=\bar\mau-  W$, then we have
 \begin{eqnarray}
 \begin{cases}
 \va\Delta\hat \mau +\lambda\va\nabla\mathrm{div}\hat\mau=\mathbf{g}-\bar {\mathbf{f}}\in L^p(Q_2(0))\quad\text{in}~Q_{2 }(0),\\
\hat \mau=0\qquad \text{on}\; T_1\cup T_2\\
\hat \mau=-W\qquad \text{on}\; \partial Q_2(0)\setminus \{T_1\cup T_2\}.\end{cases}\end{eqnarray}
By Theorem 3.8.1 in \cite{K-M-R}, for $P^*$ defined in (\ref{p0}), if $2<p<P^*$, then we have using Lemma \ref{lw} that
\begin{eqnarray}
&&\va\|\hat \mau\|_{2,p;Q_{\f32}(0)}\leq C[\|\mathbf{g}\|_{p;Q_{2}(0)}+\|\bar{\mathbf{f}}\|_{p;Q_{2}(0)}+\va\|W\|_{1,p;Q_{2}(0)}]\nonumber\\
\leq && C[\|\mathbf{f}\|_{p;\Omega}+\|\mau\|_{1,p;\Omega}+\|\tilde\rho\|_{p;\Omega}+\|\partial_{x_1}\bar\rho\|_{p;Q_{2}(0)}].\label{h3}
\end{eqnarray}
If we define $\Omega_1=(0,\f23L)\times(0,\f43)$, then combing (\ref{h4}), Lemma \ref{lw} and (\ref{h3}) we obtain that
\begin{eqnarray}
&&\va\|\partial_1\mau\|_{1,p;\Omega_1}=\va\|\partial_1\bar\mau\|_{1,p;\Omega_1}\leq \va\|\partial_1\hat\mau\|_{1,p;\Omega_1}+\va\|\partial_1W\|_{1,p;\Omega_1}\nonumber\\
\lesssim&&\|\mathbf{f}\|_{p;\Omega}+\|\mau\|_{1,p;\Omega}+\|\varrho\|_{1,p;\Omega}+\| \rho\|_{p;\Omega}+\|\partial_{x_1} \rho\|_{p;\Omega}\nonumber\\
\lesssim&&\|\mathbf{f}\|_{p;\Omega}+\|f_0\|_{p;\Omega}+\|\mau\|_{1,p;\Omega}+\|\varrho\|_{1,p;\Omega},\label{o1}
\end{eqnarray}
here we have used the mass equation $(\ref{app})_1$,  i.e.:
$$\|\rho\|_{p;\Omega}\leq \|\rho_{x_1}\|_{p;\Omega}\lesssim \|f_0\|_{p;\Omega}+\|\dv\mau\|_{p;\Omega}.$$
Similarly, if we define
$$\Omega_2=(0,\f23L)\times(\f23,2),\ \Omega_3=(\f13L,L)\times(\f23,2),\ \Omega_4=(\f13L,L)\times(0,\f43),$$
then $\Omega=\cup_{i=1}^4\Omega_i$, and estimate (\ref{o1}) also holds in $\Omega_i,i=2,3,4.$  Consequently we have
 \begin{equation}\va\|\partial_1\mau\|_{1,p;\Omega}\lesssim\|\mathbf{f}\|_{p;\Omega}+\|f_0\|_{p;\Omega}+\|\mau\|_{1,p;\Omega}+\|\varrho\|_{1,p;\Omega}.\label{du1}\end{equation}
Combining (\ref{du}), (\ref{du1}) and Lemma \ref{lp} we have
\begin{eqnarray}
&&\va\|\partial_1\mau\|_{1,p;\Omega}\lesssim\|\mathbf{f}\|_{p;\Omega}+\|f_0\|_{p;\Omega}+\|\mau\|_{1,p;\Omega}+\|\varrho\|_{1,p;\Omega}\nonumber\\
\lesssim&&\|\mathbf{f}\|_{p;\Omega}+\|f_0\|_{p;\Omega}+\va^{1+\f\sigma{10}}(\|\mau\|_{2,p;\Omega}+\|f_0\|_{W^{1,p}})+\va^{-\f32+\f2p-\f12\sigma}(\|\mathbf{f}\|+\|f_0\|)\nonumber\\
\lesssim&&\|\mathbf{f}\|_{p;\Omega}+\va^{1+\f\sigma{10}}(\|\mau\|_{2,p;\Omega}+\|f_0\|_{w^{1,p}})+\va^{-\f32+\f2p-\f12\sigma}(\|\mathbf{f}\|+\|f_0\|).\label{u1p}
\end{eqnarray}

Finally, from equation $(\ref{app})_2$, $(\ref{app})_3$, (\ref{u1p}) and Lemma \ref{rhoy} we have
\begin{eqnarray}
&&\va\|\mau_{x_2x_2}\|_{p,\Omega}\leq \|\mathbf{f}\|_{p;\Omega}+\|\mau\|_{1,p;\Omega}+\|\rho_{x_2}\|_{p;\Omega}+\va\|\partial_1\mau\|_{1,p;\Omega}\nonumber\\
\leq&&\|\mathbf{f}\|_{p;\Omega}+\va^{1+\f\sigma{10}}\|\mau\|_{2,p;\Omega}+\va^{-\f32+\f2p-\f12\sigma}(\|\mathbf{f}\|+\|f_0\|)+\va\|f_0\|_{W^{1,p}}.
\end{eqnarray}
Combing above, we finish the proof of Theorem \ref{2p}.

\subsection{Solutions to the Linear System}
As the estimates in Theorem \ref{thapp} and Theorem \ref{2p} are independent of $\delta$, by taking the limit with $\delta\rightarrow0$, we have the following theorem:
\begin{theorem}\label{thlinear}
For given $f_0\in L^2(\Omega),\ \mathbf{f}=(f_1,f_2)\in (H^{-1}(\Omega))^2$, , there exists a unique weak solution $(\rho,u,v)\in L^2(\Omega)\times H_0^1(\Omega)\times H_0^1(\Omega)$ to system (\ref{app})-(\ref{boundary.app}) with the following estimates:
\begin{equation*}\|\mau\|^2+\|\rho\|^2+|\rho(L,\cdot)|^2+\va\|\nabla\mau\|^2\leq C(\va)(\|f_0\|^2+\|\mathbf{f}\|_{H^{-1}}).\end{equation*}
Moreover, if $f_0\in W^{1,p}(\Omega),\ \mathbf{f}=(f_1,f_2)\in (L^p(\Omega))^2$, then we have $(\rho,u,v)\in W^{1,p}(\Omega)\times (W^{2,p}(\Omega))^2$
and
\begin{eqnarray*}
&&\|\mau\|^2+\|\rho\|^2+|\rho(L,\cdot)|^2+\va\|\nabla\mau\|^2\leq C_1(\|f_0\|^2+\|\mathbf{f}\|_{L^2}),\\
&&\|\rho\|_{1,p;\Omega}+\va\|\mathbf{u}\|_{2,p;\Omega}\leq
C_2[\|\mathbf{f}\|_{p;\Omega}+\va^{-\f32+\f2p-\f12\sigma}(\|\mathbf{f}\|+\|f_0\|)+\va\|f_0\|_{W^{1,p}}],
\end{eqnarray*}
where $\sigma>0$ is any constant small enough and the constants $C_1,C_2$ are independent of
$\va$.
\end{theorem}
\renewcommand{\theequation}{\thesection.\arabic{equation}}
\setcounter{equation}{0}
\section{Solutions to the Nonlinear System}

In this section we will prove the existence of solutions to the nonlinear system.
First the following Lemma gives the uniform bound of $\{(\mathbf{u}^{n},\rho^{n})\}$
in both Hilbert space and $L^p$ space.

\begin{lemma}\label{lemma3.1}
Let $\{(\mathbf{u}^{n},\rho^{n})\}$ be the sequence of solutions to
the system (\ref{n0}) with $(\mathbf{u}^{0},\rho^{0})=(0,0,0)$, then there exist  constants $M_1,M_2>0$ such that
for any $n\in Z^+$, we have
\begin{equation}\|\mau^n\|+\|\rho^n\|+|\rho^n(L,\cdot)|+\va^{\f12}\|\nabla\mau^n\|\leq M_1\va^{\f52-\f2p+\sigma},\label{n2}\end{equation}
and
\begin{eqnarray}
&&\va\|\mathbf{u}^n\|_{W^{2,p}}+\|\rho^n\|_{W^{1,p}}\leq M_2\va^{1+\f\sigma2}\label{n1}.
\end{eqnarray}

\end{lemma}

\begin{proof}
First of all, if we take
$$f_0(\rho^n,u^n,v^n)=g_0(\mau^n,\rho^n)+\bar g_0(v^n),\ \mathbf{f}(\rho^n,u^n,v^n)=\mathbf{g}(\rho^n,\mau^n)+\bar {\mathbf{g}},$$
then by Theorem \ref{thlinear}  we have
\begin{eqnarray}
&&\|\mau^{n+1}\|+\|\rho^{n+1}\|+|\rho^{n+1}(L,\cdot)|+\va^{\f12}\|\nabla\mau^{n+1}\|\nonumber\\
\leq&& C_1(\|\mathbf{f}(\rho^n,u^n,v^n)\|+\|f_0(\rho^n,u^n,v^n)\|)
\end{eqnarray}
and
\begin{eqnarray}
&&\va\|\mathbf{u}^{n+1}\|_{W^{2,p}}+\|\rho^{n+1}\|_{W^{1,p}}\nonumber\\
\leq&& C_2[\va\|f_0(\rho^n,u^n,v^n)\|_{W^{1,p}}+\va^{\f2p-\f32-\f12\sigma}(\|\mathbf{f}(\rho^n,u^n,v^n)\|+\|f_0(\rho^n,u^n,v^n)\|)\nonumber\\
&&+\|\mathbf{f}(\rho^n,u^n,v^n)\|_{L^p}].
\end{eqnarray}
where the constant $C_1,C_2$ depends only on $\Omega$ and $ p$.

Now we prove (\ref{n2}) and (\ref{n1}) by an induction argument.

 When $n=1$, by the assumption that $(\mathbf{u}^0,\rho^0)=(0,0)$ we have
 $$f_0(0,0,0)=g_{0s}+\bar {g}_{0s},\ f_1(0,0,0)=g_{1s}+\bar g_{1},\ f_2(0,0,0)=g_{2s}+\bar g_{2},$$
here $\bar {g}_{0s},\bar{\mathbf{g}}=(\bar g_1,\bar g_2)$ are defined in (\ref{linear.4})-(\ref{linear.5}). Then we have
\begin{eqnarray*}
&&\|\mau^1\|+\|\rho^1\|+|\rho^1(L,\cdot)|+\va^{\f12}\|\nabla\mau^1\|
\leq C_1(\|f_0(0,0,0)\|_{L^2}+\|\mathbf{f}(0,0,0)\|_{L^2})\\
=&& C_1(\|g_{0s}\|_{L^2}+\|\mathbf{g}_{s}\|_{L^2}+\|\bar {g}_{0s}\|_{L^2}+\|\bar{\mathbf{g}}\|_{L^2})\\
\leq&&\bar C_1(\va^2+\va^{\f74}+\va^{\f52-\f2p+ \sigma})\leq 2\bar C_1\va^{\f52-\f2p+ \sigma},
\end{eqnarray*}
and
\begin{eqnarray*}
&&\|\rho^1\|_{1,p}+\va\|\mau^1\|_{2,p} \\
\leq&&C_2[\va\|f_0(0,0,0)\|_{1,p}+\va^{\f2p-\f32-\f12\sigma}(\|f_0(0,0,0)\|_{L^2}+\|\mathbf{f}(0,0,0)\|_{L^2})+\|\mathbf{f}(0,0,0)\|_{L^p}]\\
\leq&&C_2[\va(\|g_{0s}\|_{1,p}+\|\bar g_{0s}\|_{1,p})+\va^{\f2p-\f32-\f12\sigma}(\|g_{0s}\|_{L^2}+\|\bar g_{0s}\|_{L^2}+\|\mathbf{g}_s\|_{L^2}+\|\bar{\mathbf{g}}\|_{L^2})\\
&&+\|\mathbf{g}_s\|_{L^p}+\|\bar{\mathbf{g}}_s\|_{L^p}]\\
\leq&&2\bar C_2\va^{1+\f12\sigma}.
\end{eqnarray*}
If we take $\bar C_1\leq \f {M_1}2$, $\bar C_2\leq \f {M_2}2$, then (\ref{n1}) holds for $n=1$.

 Now, assuming that for any $1\leq k\leq n$,
\begin{eqnarray*}
A_k\triangleq\|\mau^k\|+\|\rho^k\|+|\rho^k(L,\cdot)|+\va^{\f12}\|\nabla\mau^k\|\leq M_1\va^{\f52-\f2p+\sigma}
\end{eqnarray*}
and
\begin{eqnarray*}
B_k\triangleq\|\rho^k\|_{1,p}+\va\|\mau^k\|_{2,p}\leq M_2\va^{1+\f12\sigma}.
\end{eqnarray*}
Then direct computation shows that
\begin{eqnarray*}
&&\|\mathbf{g}_{r}(\rho^n,u^n,v^n)\|+\|g_{0r}(\rho^n,u^n,v^n)\|\\
\leq&&C[\va\|\mau^n\|_{H^1}+\va\|\rho^n\|+\va^{\f12}\|v^n\|+\|\mau^n\|_{L^\infty}(\|\nabla\mau^n\|+\|\rho^n\|)\nonumber\\
 &&+\|\rho^n\|_{L^\infty}\|\mau^n\|_{H^1}+\|\rho^n\|_{L^\infty}\|\mau^n\|_{L^\infty}\|\nabla\mau^n\|+\|\rho^n\|_{L^\infty}\|\nabla\rho^n\|]\nonumber\\
  \leq&&C[\va^{\f12}A_n+\va B^{2,n}+\va^{-1-2\sigma}A^2_n+\va^{-\f12-2\sigma}B^{2,n}A_n+\va^{-2\sigma}B^{2,n}(\va^{-\f12}A_n+\va^{-1-\sigma}A_n^2)\nonumber\\
  &&+\va^{-\sigma}B^{2,n}B^{2,n}]\nonumber\\
  \leq&&C\va^{2-10\sigma}.
\end{eqnarray*}
Consequently we have
\begin{eqnarray}
&&\|\mau^{n+1}\|+\|\rho^{n+1}\|+|\rho^{n+1}(L,\cdot)|+\va^{\f12}\|\nabla\mau^{n+1}\|\nonumber\\
\leq&& C_2[\va\|f_0(\rho^n,u^n,v^n)\|_{W^{1,p}}+\va^{\f2p-\f32-\f12\sigma}(\|\mathbf{f}(\rho^n,u^n,v^n)\|+\|f_0(\rho^n,u^n,v^n)\|)\nonumber\\
&&+\|\mathbf{f}(\rho^n,u^n,v^n)\|_{L^p}]\nonumber\\
\leq&&C_2[\va(\|g_{0s}\|_{1,p}+\|\bar g_{0s}\|_{1,p})+\va^{\f2p-\f32-\f12\sigma}(\|g_{0s}\|_{L^2}+\|\bar g_{0s}\|_{L^2}+\|\mathbf{g}_s\|_{L^2}+\|\bar{\mathbf{g}}\|_{L^2})\nonumber\\
&&+\|\mathbf{g}_s\|_{L^p}+\|\bar{\mathbf{g}}_s\|_{L^p}]+\|\mathbf{g}_{r}(\rho^n,u^n,v^n)\|+\|g_{0r}(\rho^n,u^n,v^n)\|)\nonumber\\
 \leq&& C_2\va^{\f52-\f2p+ \sigma}+C_2\va^{2-2\sigma} \leq M_1\va^{\f52-\f2p+ \sigma},
\end{eqnarray}
and
\begin{eqnarray*}
&&\|\rho^{n+1}\|_{1,p}+\va\|\mau^{n+1}\|_{2,p}\\
\leq&& C_2[\va\|f_{0}(\rho^n,u^n,v^n)\|_{1,p}+\va^{\f2p-\f32-\f12\sigma}(\|f_{0}(\rho^n,u^n,v^n)\|_{L^2}+\|\mathbf{f}(\rho^n,u^n,v^n)\|_{L^2})\\
&&+\|\mathbf{f}(\rho^n,u^n,v^n)\|_{L^p}]\\
\leq&&C_2[\va\|g_{0s}\|_{1,p}+\va^{\f2p-\f32-\f12\sigma}(\|g_{0s}\|_{L^2}+\|\mathbf{g}_s\|_{L^2})+\|\mathbf{g}(\rho^n,u^n,v^n)\|_{L^p}\\
&&+\va\|g_{0r}(\rho^n,u^n,v^n)\|_{1,p}+\va^{\f2p-\f32-\f12\sigma}(\|g_{0r}(\rho^n,u^n,v^n)\|_{L^2}+\|\mathbf{g}_r(\rho^n,u^n,v^n)\|_{L^2})\\
\leq&&C_2[\va^{1+\f12\sigma}+\va^2(\|\rho^n\|_{1,p}+\|\mau^n\|_{2,p})+\va(\|\rho^n\|_{L^\infty}\|\nabla^2\mau^n\|_{L^p}+\|\nabla\rho^n\|_{L^p}\|\nabla\mau^n\|_{L^\infty})\\
&&+\va^2(\|\rho^n\|_p+\|\mau^n\|_{1,p})+\va(\|\mau^n\|_{L^\infty}+\|\rho^n\|_{L^\infty})(\|\nabla\rho^n\|_{L^p}+\|\mau^n\|_{1,p})\\
&&+\va\|\rho^n\|_{L^\infty}\|\mau^n\|_{L^\infty}\|\mau^n\|_{1,p}]\\
\leq&&C_2[\va^{1+\f12\sigma}+\va^2B_n+B_n^2+B_n^2A_n+A_n^2B_n]\\
\leq&&M_2\va^{1+\f12\sigma}.
\end{eqnarray*}
Thus we have finished the proof of the Lemma.

\end{proof}
If we denote by
$$\mathcal{X}=\{u\in H^1,|\|u\|_{\mathcal{X}}<\infty\},\ \ \|u\|_{\mathcal{X}}=\va^{\f12}\|\nabla u\|_{L^2}+\|u\|_{L^2},$$
then based on Lemma \ref{lemma3.1}, we can prove that $\{(\mathbf{u}^{n},\rho^{n})\}$ is a Cauchy sequence in $(\mathcal{X})^2\times L^2$. More precisely, we have the following Lemma:
\begin{lemma}\label{cs}
Let $\{(\mathbf{u}^{n},\rho^{n})\}_{n=1}^\infty$ be the sequence of solutions to system (\ref{n0}) with initial data $(\mathbf{u}^{0},\rho^{0})=(0,0,0)$, then we have
\begin{eqnarray}
&&\va^{\f12}\|\nabla(\mathbf{u}^{n+1}-\mathbf{u}^{n})\|_{L^2}+\|\mathbf{u}^{n+1}-\mathbf{u}^{n}\|_{L^2}+\|\rho^{n+1}-\rho^{n}\|_{L^2}\nonumber\\
\leq&&\f12 \Big(\va^{\f12} \|\nabla(\mathbf{u}^{n}-\mathbf{u}^{n-1})\|_{L^2}+ \|\mathbf{u}^{n}-\mathbf{u}^{n-1}\|_{L^2}+\|\rho^{n}-\rho^{n-1}\|_{L^2}\Big).
\end{eqnarray}
\end{lemma}

\begin{proof} A straightforward calculation gives
\begin{eqnarray*}
&& \dv(\mau^{n+1}-\mau^n)+(\mau_s+\mau^{n})\cdot\nabla(\rho^{n+1}-\rho^n)\nonumber\\
=&&(v^n-v^{n-1})(\rho_{0}'-\rho_{sx_2}(0,x_2))+g_{0r}(\mathbf{u}^{n},\rho^{n})-g_{0r}(\mathbf{u}^{n-1},\rho^{n-1})\\
&&- (\mathbf{u}^{n}-\mathbf{u}^{n-1})\cdot\nabla\rho^{n}\nonumber\\
\triangleq&&\hat g_0,\nonumber\\
&&\  u_s(u^{n+1}-u^n)_{x_1}+u_{sx_2}(v^{n+1}-v^n)-\va\Delta (u^{n+1}-u^n)
-\va\partial_{x_1}\mathrm{div}(\mau^{n+1}-\mau^n)\nonumber\\
&&\quad+c^2(\rho^{n+1}-\rho^n)_{x_1}\nonumber\\
=&&g_{1r}(\mathbf{u}^{n},\rho^{n})-g_{1r}(\mathbf{u}^{n-1},\rho^{n-1})\nonumber\\
\triangleq&& \hat g_1,\nonumber\\
&&\  u_s(v^{n+1}-v^n)_{x_1}-\va\Delta (v^{n+1}-v^n)
-\va\partial_{x_2}\mathrm{div}(\mau^{n+1}-\mau^n)+c^2(\rho^{n+1}-\rho^n)_{x_2}\nonumber\\
=&&g_{2r}(\mathbf{u}^{n},\rho^{n})-g_{2r}(\mathbf{u}^{n-1},\rho^{n-1})\\
\triangleq&&\hat g_2,
\end{eqnarray*}
with boundary condition
$$\rho^{n+1}-\rho^{n}|_{x_1=0}=0,\
\mau^{n+1}-\mau^n|_{\partial\Omega}=0.$$
Take $t=1$ in (\ref{L2}) we have
\begin{eqnarray*}
&&\|\mau^{n+1}-\mau^n\|^2+ \|\rho^{n+1}-\rho^n\|^2+\va\|\nabla(\mau^{n+1}-\mau^n)\|^2\nonumber\\
\leq&& C(\|\hat g_0\|^2+| ((L-x)(\mau^{n+1}-\mau^n),\hat{\mathbf{g}})|+|(\mau^{n+1}-\mau^n,\hat{\mathbf{g}})|),\\
\leq&& C(\|\hat g_0\|^2+| ((L-x)(v^{n+1}-v^n),f^{2,n})|+\|\mau^{n+1}-\mau^n\|(\|\hat g_1\|+\|\hat g_2-f^{2,n}\|))
\end{eqnarray*}
here $\hat{\mathbf{g}}=(\hat g_1,\hat g_2)$, $f^{2,n}=[c^2-p'(\rho^{\va,n})]\rho^n_{x_2}-[c^2-p'(\rho^{\va,n-1})]\rho^{n-1}_{x_2}$
and $\rho^{\va,n}=\rho_s+\rho^n$. Then by  direct computation we have for $v^n-v^{n-1}\in H^1_0\cap W^{2,p}(\Omega)$,
\begin{eqnarray*}
&&|(f^{2,n},v^n-v^{n-1})|\\
=&&|\int_\Omega\{[c^2-p'(\rho^{\va,n})]\rho^n_{x_2}-[c^2-p'(\rho^{\va,n-1})]\rho^{n-1}_{x_2}\}(v^n-v^{n-1})dx_1dx_2|\\
\leq&&\int_\Omega[|p'(\rho^{\va,n})-p'(\rho^{\va,n-1})||\rho^n_{x_2}||v^n-v^{n-1}|dx_1dx_2\\
&&+|\int_\Omega[c^2-p'(\rho^{\va,n-1})](\rho^n-\rho^{n-1})_{x_2}(v^n-v^{n-1})dx_1dx_2|\\
\leq&&\int_\Omega[|p'(\rho^{\va,n})-p'(\rho^{\va,n-1})||\rho^n_{x_2}||v^n-v^{n-1}|dx_1dx_2\\
&&+|\int_\Omega[(c^2-p'(\rho^{\va,n-1}))(v^n-v^{n-1})]_{x_2}(\rho^n-\rho^{n-1})dx_1dx_2|\\
\leq&&\|(\rho^n-\rho^{n-1})\|_{L^2}[\|\rho_{x_2}^{n-1}\|_{L^p}\|v^n-v^{n-1}\|_{L^{\f{2p}{p-2}}}+\|\rho^{n-1}\|_{L^\infty}\|(v^n-v^{n-1})_{x_2}\|_{L^2}]\\
\leq&&\va\|(\rho^n-\rho^{n-1})\|_{L^2}\|v^n-v^{n-1}\|_{H^1}\\
\leq&&\va^{\f12}(\|(\rho^n-\rho^{n-1})\|_{L^2}^2+\va\|\nabla(v^n-v^{n-1})\|_{L^2}^2+\|v^n-v^{n-1}\|_{L^2}^2)
\end{eqnarray*}
 Besides, we have
\begin{eqnarray*}
&&\|\hat g_0\|_{L^2}+|\hat g_{1}|_{L^2}+|\hat g_{2}-f^{2,n}|_{L^2}\\
\lesssim&&\|\rho^n\dv\mau^n-\rho^{n-1}\dv\mau^{n-1}\|_{L^2}+\|(\mathbf{u}^{n}-\mathbf{u}^{n-1})\cdot\nabla\rho^{n}\|_{L^2}\\
&&+\va\|\rho^n-\rho^{n-1}\|_{L^2}+(\|\rho^{n-1}\|_{L^\infty}+\|\mau^{n-1}\|_{L^\infty})\|\mau^n-\mau^{n-1}\|_{H^1}\\
&&+\va\|\mathbf{u}^{n}-\mathbf{u}^{n-1}\|_{H^1}+(\|\mau^n-\mau^{n-1}\|_{L^2}+\|\rho^n-\rho^{n-1}\|_{L^2})\|\nabla\mau^n\|_{L^\infty}\\
&&+\|\rho^n\mau^n\cdot\nabla\mau^n-\rho^{n-1}\mau^{n-1}\cdot\nabla\mau^{n-1}\|_{L^2}\\
\leq&& \va^{\f\sigma2} (\va^{\f12}\|\nabla(\mathbf{u}^{n+1}-\mathbf{u}^{n})\|_{L^2}+\|\rho^n-\rho^{n-1}\|_{L^2}+\|\mathbf{u}^{n+1}-\mathbf{u}^{n}\|_{L^2}),
\end{eqnarray*}
 Combing above,  we finish  the proof of Theorem \ref{app}.
\end{proof}

\textbf{Proof of Theorem \ref{main}.} \ Lemma \ref{cs} implies that
$\{(\mathbf{u}^{n},\rho^{n})\}$ is a Cauchy sequence in
$(\mathcal{X})^2\times L^2$. Hence, with the help of Lemma \ref{lemma3.1} we see that
there exists a unique $({\mathbf{u}},\rho)\in W^{2,p}\times W^{1,p}$,
 such that, such that for any $1<p'<p$,
$$\mathbf{u}^{n}\rightarrow\mathbf{u}\;\quad\text{strongly in }\; W^{2,p'}(\Omega), $$
$$\rho^{n}\rightarrow\rho\;\quad\text{strongly in }\; W^{1,p'}(\Omega).$$

 Taking limit on both sides of (\ref{n0})  we conclude that
 $(\mathbf{u},\rho)$ satisfies the nonlinear system (\ref{0.3})-(\ref{0.5})  with boundary condition (\ref{b}). Moreover, the estimate (\ref{th.0.1})-(\ref{th.0.4})  follows  immediately from (\ref{n2})-(\ref{n1}) and  Sobolev imbedding Theorem.

Next, if $(\mathbf{u},\rho)$ and $(\hat{\mathbf{u}},\hat{\rho})$ are
two solutions of the system  (\ref{0.3})-(\ref{0.5}), then by a process similar to that used
in Lemma \ref{cs}, we infer that
\begin{eqnarray*}
&&\va^{\f12}\|\nabla(\mathbf{u}-\hat{\mathbf{u}})\|_{2;\Omega}+\|\mathbf{u}-\hat{\mathbf{u}}\|_{2;\Omega}+\|\rho-\hat{\rho}\|_{2;\Omega}\\
\leq&&\f12(\va^{\f12}\|\nabla(\mathbf{u}-\hat{\mathbf{u}})\|_{2;\Omega}+\|\mathbf{u}-\hat{\mathbf{u}}\|_{2;\Omega}+\|\rho-\hat{\rho}\|_{2;\Omega}).
\end{eqnarray*}
Consequently, we obtain the uniqueness of the solutions. This completes the proof of Theorem \ref{main}.

\noindent\textbf{Acknowledgements.} Chunhui Zhou is supported by  the Jiangsu Provincial Scientific Research Center of Applied Mathematics under Grant No.BK20233002.

\renewcommand{\theequation}{\thesection.\arabic{equation}}
\setcounter{equation}{0}
\appendix

\section{Weight estimates for a linear parabolic equation}
For simplicity, we first consider the following parabolic system:
\begin{eqnarray} \label{pa1}
\begin{cases}
A_0 \p_{x_1} u - \p_{YY}u = 0, \ 0<x_1<L,\ 0<Y<+\infty,\\
u|_{x_1 = 0} = u^{0,Y}(Y) ,\ 0<Y<+\infty,\\
 u|_{Y = 0} = u_0(x_1), u|_{Y \rightarrow +\infty} = 0, \ 0<x_1<L,
\end{cases}
\end{eqnarray}
here $A_0>0$ is a given constant. We define a cutoff function $w(Y)\in C^\infty(0,\infty)$, $0\leq w(Y)\leq 1$, $\|w\|_{C^5}\leq C$, and
\begin{equation}w(Y)=\begin{cases}1,\ Y>4,\\0,\ 0<Y<3.\end{cases}\label{A1}\end{equation}
We also denote by $\Omega_p=(0,L)\times(0,\infty)$.

\begin{lemma} \label{ba}Assume that $M\in \mathbb{Z}^+$, $u_0(x)\in W^{M,p}(0,L)$,  $u^{0,Y}(Y)$ is a smooth function in $[0,\infty)$ with compact support in $[0,2]$,
and the  compatibility conditions on $(0,0)$ hold:
$$\partial_{x_1}^k u_{0}(0)=\f1{A_0^k}\partial_Y^{2k}u^{0,Y}(0),\ k=0,...,M-1,$$
then there exists a unique solution $u$ to system (\ref{pa1})
with the following estimate:
\begin{eqnarray}
&&\|(1+Y)^m w(Y)\nabla ^ju \|_{L^\infty}+\|\partial_{x_1}^k\partial_Y^{l}u\|_{p,\Omega_p}\nonumber\\
 \leq&& C(m,j,k)(|u_0|_{M,p}+\|u^{0,Y}\|_{2M,p}),\  \text{ for }\  0 \le 2k+ l \le 2M, j\in N\label{ba.1}.
\end{eqnarray}
where the constant $C$ does not depend on $Y$.

\end{lemma}

\begin{proof}  First of all, by the $L^p$ theory of linear parabolic equations, there exists a  unique solution $u$  to system (\ref{pa1}) with $\partial_x^k\partial_Y^{l}u\in L^p(\Omega_p)$, $0\leq 2k+l\leq2M$, and   the following estimate holds:
$$\|\partial_{x_1}^k\partial_Y^{l}u\|_{p,\Omega_p}\lesssim |u_0|_{M,p}+\|u^{0,Y}\|_{2M,p},$$
and
$$\|w(Y)u\|_{j,p;\Omega_p}\leq C(j)(|u_0|_{M,p}+\|u^{0,Y}\|_{2M,p}),\ \text{for any}\ j\in N.$$

 To get weighted  estimates, we  multiply equation (\ref{pa1}) with $ u_{x_1} w^2(Y)(1+Y)^{2m}$ to have
 \begin{eqnarray*}
 &&A_0\int_0^\infty w^2(Y)(1+Y)^{2m} u_{x_1}^2dy+\f12\f{d}{dx_1}\int_0^\infty w^2(Y)(1+Y)^{2m}( u_{Y})^2dy\nonumber\\
 =&&-\int_0^\infty [2mw^2(Y)(1+Y)^{2m-1}+2w(Y)w'(Y)(1+Y)^{2m}] u_{x_1} u_{Y}dy\nonumber\\
\leq &&2m\|w(Y)(1+Y)^{m-1} u_{x_1}(x_1,\cdot)\|_{L^2_Y}\|w(Y)(1+Y)^mu_Y(x_1,\cdot)\|_{L^2_Y}\nonumber\\
&&+C\|u_{x_1}\|_{L^2_Y}\|u_Y\|_{L^2_Y}.\nonumber
 \end{eqnarray*}
Using the fact that $L\ll1$, we have
 \begin{eqnarray*}
&& \|w(Y)(1+Y)^m u_{x_1}\|_{L^2(\Omega)}+\sup_{x_1\in[0,L]}\|w(Y)(1+Y)^m u_Y(x_1,\cdot)\|_{L^2_y}\nonumber\\
\leq &&C(m)\|w(Y)(1+Y)^{m-1} u_{x_1}(x_1,\cdot)\|_{L^2}^2+C\|u_{x_1}\|_{L^2_Y([2,5])}\|u_Y\|_{L^2_Y([2,5])}.
 \end{eqnarray*}
 By a process of induction, we have
  \begin{eqnarray*}
&& \|w(Y)(1+Y)^m u_{x_1}\|_{L^2(\Omega)}+\sup_{x_1\in[0,L]}\|w(Y)(1+Y)^m u_Y(x_1,\cdot)\|_{L^2_y}\nonumber\\
\leq &&C(m)(\|w(Y) u_{x_1}\|_{L^2}^2+\|w(Y)u_Y\|_{L^2}^2)+C\|u_{x_1}\|_{L^2_Y([2,5])}\|u_Y\|_{L^2_Y([2,5])}\nonumber\\
\lesssim&& |u_0|_{M,p}+\|u^{0,Y}\|_{2M,p}.
 \end{eqnarray*}
 Similarly, by taking derivative of equation (\ref{pa1}) with respect to $x_1$ and repeat the process above we  obtain for any $m,j\in N$,
  \begin{eqnarray*}
 \|w(Y)(1+Y)^m\nabla^{j+2} u\|_{L^2(\Omega)}\leq C(m,j)|u_0|_{M,p}+\|u^{0,Y}\|_{2M,p}.
 \end{eqnarray*}
 The Sobolev imbedding inequality implies that
\begin{eqnarray*}\|w(Y)(1+Y)^m\nabla^j u\|_{L^\infty(\Omega)}\leq C(m,j)|u_0|_{M,p}+\|u^{0,Y}\|_{2M,p}.\end{eqnarray*}
Thus we have finished the proof.
\end{proof}

\end{document}